\documentclass[12pt]{article}
\usepackage{latexsym,amsfonts,amsmath,amssymb}
\newtheorem{theorem}{Theorem}
\newtheorem{corollary}[theorem]{Corollary}
\newtheorem{sublemma}{Lemma}[theorem]
\newtheorem{lemma}[theorem]{Lemma}
\newtheorem{question}[theorem]{Question}
\newtheorem{observation}[theorem]{Observation}
\newtheorem{claim}[theorem]{Claim}
\newtheorem{subclaim}{Claim}[sublemma]
\newtheorem{conjecture}[theorem]{Conjecture}
\newtheorem{fact}[theorem]{Fact}
\newtheorem{definition}[theorem]{Definition}
\newtheorem{remark}[theorem]{Remark}
\newtheorem{example}[theorem]{Example}
\newtheorem{exercise}{Exercise}[section]
\def\Theorem #1.#2 #3\par{\setbox1=\hbox{#1}\ifdim\wd1=0pt
   \begin{theorem}{\rm #2} #3\end{theorem}\else
   \newtheorem{#1}[theorem]{#1}\begin{#1}\label{#1}{\rm #2} #3\end{#1}\fi}
\def\Corollary #1.#2 #3\par{\setbox1=\hbox{#1}\ifdim\wd1=0pt
   \begin{corollary}{\rm #2} #3\end{corollary}\else
   \newtheorem{#1}[theorem]{#1}\begin{#1}\label{#1}{\rm #2} #3\end{#1}\fi}
\def\Lemma #1.#2 #3\par{\setbox1=\hbox{#1}\ifdim\wd1=0pt
   \begin{lemma}{\rm #2} #3\end{lemma}\else
   \newtheorem{#1}[theorem]{#1}\begin{#1}\label{#1}{\rm #2} #3\end{#1}\fi}
\def\SubLemma #1.#2 #3\par{\setbox1=\hbox{#1}\ifdim\wd1=0pt
   \begin{sublemma}{\rm #2} #3\end{sublemma}\else
   \newtheorem{#1}{#1}[theorem]\begin{#1}\label{#1}{\rm #2} #3\end{#1}\fi}
\def\Question #1.#2 #3\par{\setbox1=\hbox{#1}\ifdim\wd1=0pt
   \begin{question}{\rm #2} #3\end{question}\else
   \newtheorem{#1}[theorem]{#1}\begin{#1}\label{#1}{\rm #2} #3\end{#1}\fi}
\def\Observation #1.#2 #3\par{\setbox1=\hbox{#1}\ifdim\wd1=0pt
   \begin{observation}{\rm #2} #3\end{observation}\else
   \newtheorem{#1}[theorem]{#1}\begin{#1}\label{#1}{\rm #2} #3\end{#1}\fi}
\def\Claim #1.#2 #3\par{\setbox1=\hbox{#1}\ifdim\wd1=0pt
   \begin{claim}{\rm #2} #3\end{claim}\else
   \newtheorem{#1}[theorem]{#1}\begin{#1}\label{#1}{\rm #2} #3\end{#1}\fi}
\def\SubClaim #1.#2 #3\par{\setbox1=\hbox{#1}\ifdim\wd1=0pt
   \begin{subclaim}{\rm #2} #3\end{subclaim}\else
   \newtheorem{#1}{#1}[sublemma]\begin{#1}\label{#1}{\rm #2} #3\end{#1}\fi}
\def\Conjecture #1.#2 #3\par{\setbox1=\hbox{#1}\ifdim\wd1=0pt
   \begin{conjecture}{\rm #2} #3\end{conjecture}\else
   \newtheorem{#1}[theorem]{#1}\begin{#1}\label{#1}{\rm #2} #3\end{#1}\fi}
\def\Fact #1.#2 #3\par{\setbox1=\hbox{#1}\ifdim\wd1=0pt
   \begin{fact}{\rm #2} #3\end{fact}\else
   \newtheorem{#1}[theorem]{#1}\begin{#1}\label{#1}{\rm #2} #3\end{#1}\fi}
\def\Definition #1.#2 #3\par{\setbox1=\hbox{#1}\ifdim\wd1=0pt
   \begin{definition}{\rm #2} {\rm #3}\end{definition}\else
   \newtheorem{#1}[theorem]{#1}\begin{#1}\label{#1}{\rm #2} {\rm #3}\end{#1}\fi}
\def\Remark #1.#2 #3\par{\setbox1=\hbox{#1}\ifdim\wd1=0pt
   \begin{remark}{\rm #2} {\rm #3}\end{remark}\else
   \newtheorem{#1}[theorem]{#1}\begin{#1}\label{#1}{\rm #2} {\rm #3}\end{#1}\fi}
\def\Example #1.#2 #3\par{\setbox1=\hbox{#1}\ifdim\wd1=0pt
   \begin{example}{\rm #2} #3\end{example}\else
   \newtheorem{#1}[theorem]{#1}\begin{#1}\label{#1}{\rm #2} #3\end{#1}\fi}
\def\Exercise #1.#2 #3\par{\setbox1=\hbox{#1}\ifdim\wd1=0pt
   {\footnotesize\begin{exercise}{\rm #2} {\rm #3}\end{exercise}}\else
   \newtheorem{#1}[section]{#1}{\footnotesize\begin{#1}\label{#1}{\rm #2} {\rm #3}\end{#1}}\fi}
\def\QuietTheorem #1.#2 #3\par{\setbox1=\hbox{#1}\ifdim\wd1=0pt\proclaim{Theorem {\rm #2}}{#3}\else\proclaim{#1 {\rm #2}}{#3}\fi}
\newcommand{\proclaim}[2]{\smallskip\noindent{\bf #1} {\sl#2}\par\smallskip}
\def\Proclaim #1.#2 #3\par{\proclaim{#1 {\rm #2}}{#3}}
\newenvironment{proof}{\noindent}{\kern2pt\QEDbox\par\bigskip}
\def\Proof#1: {\setbox1=\hbox{#1}\ifdim\wd1=0pt\begin{proof}{\bf Proof: }\else\medskip\begin{proof}{\bf #1: }\fi}
\newcommand{\QED}{\end{proof}}
\def\BF#1.{{\bf #1.}}
\def\Abstract #1\par{\begin{quotation}{\singlespaced\footnotesize{\noindent{\bf Abstract.~}#1}}\end{quotation}}
\def\Title #1\par{\title{#1}\maketitle}
\def\Author #1\par{\author{#1}}
\def\Acknowledgement#1\par{\thanks{#1}}
\def\Chapter #1\par{\chapter{#1}}
\def\Section #1\par{\section{#1}}
\def\QuietSection #1\par{\section*{#1}}
\def\SubSection #1\par{\subsection{#1}}
\def\SubSubSection #1\par{\subsubsection{#1}}
\def\MidTitle #1\par{\bigskip\goodbreak\centerline{\small\bf #1}\bigskip\noindent}

\newcommand{\singlespaced}{\baselineskip=15pt}
\def\bottomnote #1\par{{\renewcommand{\thefootnote}{}\footnotetext{#1}}}
\newcommand{\F}{{\mathbb F}}
\renewcommand{\P}{{\mathbb P}}
\newcommand{\Q}{{\mathbb Q}}

\newcommand{\Ftail}{\F_{\!\scriptscriptstyle\rm tail}}

\newcommand{\ftail}{f_{\!\scriptscriptstyle\rm tail}}

\newcommand{\Gtail}{G_{\!\scriptscriptstyle\rm tail}}

\newcommand{\Ptail}{\P_{\!\scriptscriptstyle\rm tail}}

\newcommand{\id}{\mathop{\hbox{\small id}}}
\newcommand{\one}{\mathop{1\hskip-3pt {\rm l}}}
\newfont{\msam}{msam10 at 12pt}
\newcommand{\from}{\mathbin{\vbox{\baselineskip=3pt\lineskiplimit=0pt
                         \hbox{.}\hbox{.}\hbox{.}}}}
\newcommand{\of}{\subseteq}

\newcommand{\set}[1]{\{\,{#1}\,\}}

\newcommand{\compose}{\circ}

\newcommand{\elesub}{\prec}

\newcommand{\dom}{\mathop{\rm dom}}

\newcommand{\ran}{\mathop{\rm ran}}
\newcommand{\add}{\mathop{\rm add}}

\newcommand{\Ult}{\mathop{\rm Ult}}

\newcommand{\image}{\mathbin{\hbox{\tt\char'42}}}
\newcommand{\plus}{{+}}

\newcommand{\restrict}{\upharpoonright}
\newcommand{\satisfies}{\models}
\newcommand{\forces}{\Vdash}

\newcommand{\cross}{\times}

\newcommand{\union}{\cup}

\newcommand{\intersect}{\cap}

\newcommand{\trianglelt}{\lhd}

\newcommand{\LaverDiamond}{\mathop{\hbox{\line(0,1){10}\line(1,0){8}\line(-4,5){8}}\hskip 1pt}\nolimits}
\newcommand{\LD}{\LaverDiamond}
\newcommand{\LDLD}{\mathop{\hbox{\,\line(0,1){8}\!\line(0,1){10}\line(1,0){8}\line(-4,5){8}}\hskip 1pt}\nolimits}

\newcommand{\LDwc}{\LD^{\hbox{\!\!\tiny wc}}}
\newcommand{\LDunf}{\LD^{\hbox{\!\!\tiny unf}}}
\newcommand{\LDthetaunf}{\LD^{\hbox{\!\!\tiny $\theta$-unf}}}
\newcommand{\LDsunf}{\LD^{\hbox{\!\!\tiny sunf}}}
\newcommand{\LDthetasunf}{\LD^{\hbox{\!\!\tiny $\theta$-sunf}}}
\newcommand{\LDmeas}{\LD^{\hbox{\!\!\tiny meas}}}
\newcommand{\LDstr}{\LD^{\hbox{\!\!\tiny strong}}}
\newcommand{\LDsuperstrong}{\LD^{\hbox{\!\!\tiny superstrong}}}
\newcommand{\LDram}{\LD^{\hbox{\!\!\tiny Ramsey}}}
\newcommand{\LDstrc}{\LD^{\hbox{\!\!\tiny str compact}}}
\newcommand{\LDsc}{\LD^{\hbox{\!\!\tiny sc}}}
\newcommand{\LDext}{\LD^{\hbox{\!\!\tiny ext}}}

\newcommand{\LDlambdaahuge}{\LD^{\hbox{\!\!\tiny $\lambda$-ahuge}}}
\newcommand{\LDsahuge}{\LD^{\hbox{\!\!\tiny super ahuge}}}

\newcommand{\LDlambdahuge}{\LD^{\hbox{\!\!\tiny $\lambda$-huge}}}
\newcommand{\LDshuge}{\LD^{\hbox{\!\!\tiny superhuge}}}

\newcommand{\LDlambdanhuge}{\LD^{\hbox{\!\!\tiny $\lambda$ $n$-huge}}}
\newcommand{\LDsnhuge}{\LD^{\hbox{\!\!\tiny super $n$-huge}}}
\newcommand{\LDthetasc}{\LD^{\hbox{\!\!\tiny $\theta$-sc}}}
\newcommand{\LDkappasc}{\LD^{\hbox{\!\!\tiny $\kappa$-sc}}}
\newcommand{\LDthetastr}{\LD^{\hbox{\!\!\tiny $\theta$-strong}}}
\newcommand{\LDthetastrc}{\LD^{\hbox{\!\!\tiny $\theta$-str compact}}}
\newcommand{\LDstar}{\LD^{\hbox{\!\!\tiny$\star$}}}
\newcommand{\smalllt}{\mathrel{\mathchoice{\raise2pt\hbox{$\scriptstyle<$}}{\raise1pt\hbox{$\scriptstyle<$}}{\scriptscriptstyle<}{\scriptscriptstyle<}}}
\newcommand{\smallleq}{\mathrel{\mathchoice{\raise2pt\hbox{$\scriptstyle\leq$}}{\raise1pt\hbox{$\scriptstyle\leq$}}{\scriptscriptstyle\leq}{\scriptscriptstyle\leq}}}
\newcommand{\ltkappa}{{{\smalllt}\kappa}}
\newcommand{\leqkappa}{{{\smallleq}\kappa}}
\newcommand{\leqgamma}{{{\smallleq}\gamma}}

\newcommand{\ltlambda}{{{\smalllt}\lambda}}
\newcommand{\ltomega}{{{\smalllt}\omega}}

\newcommand{\leqtheta}{{{\smallleq}\theta}}

\newcommand{\card}[1]{{|#1|}}
\def\boolval#1{\mathopen{\lbrack\!\lbrack}\,#1\,\mathclose{\rbrack\!
        \rbrack}}
\newcommand{\UnderTilde}[1]{{\setbox1=\hbox{$#1$}\baselineskip=0pt\vtop{\hbox{$#1$}\hbox to\wd1{\hfil$\sim$\hfil}}}{}}
\newcommand{\Undertilde}[1]{{\setbox1=\hbox{$#1$}\baselineskip=0pt\vtop{\hbox{$#1$}\hbox to\wd1{\hfil$\scriptstyle\sim$\hfil}}}{}}
\newcommand{\undertilde}[1]{{\setbox1=\hbox{$#1$}\baselineskip=0pt\vtop{\hbox{$#1$}\hbox to\wd1{\hfil$\scriptscriptstyle\sim$\hfil}}}{}}
\newcommand{\UnderdTilde}[1]{{\setbox1=\hbox{$#1$}\baselineskip=0pt\vtop{\hbox{$#1$}\hbox to\wd1{\hfil$\approx$\hfil}}}{}}
\newcommand{\Underdtilde}[1]{{\setbox1=\hbox{$#1$}\baselineskip=0pt\vtop{\hbox{$#1$}\hbox to\wd1{\hfil\scriptsize$\approx$\hfil}}}{}}

\newcommand{\st}{\mid}
\renewcommand{\th}{{\hbox{\scriptsize th}}}
\newcommand{\Iff}{\mathrel{\leftrightarrow}}

\newcommand{\iso}{\cong}
\def\<#1>{\langle#1\rangle}

\newcommand{\QEDbox}{\fbox{}}

\newcommand{\TC}{\mathop{\hbox{\sc tc}}}
\newcommand{\ORD}{\mathop{\hbox{\sc ord}}}

\newcommand{\ZFC}{\hbox{\sc zfc}}

\newcommand{\GCH}{\hbox{\sc gch}}

\newcommand{\AC}{\hbox{\sc ac}}

\newcommand{\Reg}{\hbox{Reg}}
\newcommand{\factordiagramup}[6]{$$\begin{array}{ccc}
#1&\raise3pt\vbox{\hbox to60pt{\hfill$\scriptstyle
#2$\hfill}\vskip-6pt\hbox{$\vector(4,0){60}$}}&#3\\ \vbox
to30pt{}&\raise22pt\vtop{\hbox{$\vector(4,-3){60}$}\vskip-22pt\hbox
to60pt{\hfill$\scriptstyle #4\qquad$\hfill}}
     &\ \ \lower22pt\hbox{$\vector(0,3){45}$}\ {\scriptstyle #5}\\
\vbox to15pt{}&&#6\\
\end{array}$$}
\newcommand{\factordiagram}[6]{$$\begin{array}{ccc}
#1&&\\ \ \ \raise22pt\hbox{$\vector(0,-3){45}$}\ {\scriptstyle #2}
&\raise22pt\hbox{$\vector(2,-1){90}$}\raise5pt\llap{$\scriptstyle#3$\qquad\quad}&\vbox
to25pt{}\\ #4&\raise3pt\vbox{\hbox to90pt{\hfill$\scriptstyle
#5$\hfill}\vskip-6pt\hbox{$\vector(4,0){90}$}}&#6\\
\end{array}$$}
\newcommand{\df}{\it} 
\begin{document}
\author{Joel David Hamkins\\
\normalsize\sc Georgia State University {\it \&}\\
\normalsize\sc The City University of New York\thanks{Specifically, The College of Staten Island of CUNY
and The CUNY Graduate Center.}
}
\date{}
\bottomnote MSC Subject Codes: 03E55, 03E35, 03E05. Keywords: Large cardinals, Laver diamond, Laver
function, fast function, diamond sequence. My research has been supported by grants from the PSC-CUNY
Research Foundation and the NSF. I would like to thank Mirna Dzamonja for helpful conversations concerning
the material in this article.

\Title A class of strong diamond principles

\Abstract In the context of large cardinals, the classical diamond principle $\Diamond_\kappa$ is easily
strengthened in natural ways. When $\kappa$ is a measurable cardinal, for example, one might ask that a
$\Diamond_\kappa$ sequence anticipate every subset of $\kappa$ not merely on a stationary set, but on a set
of normal measure one. This is equivalent to the existence of a function $\ell\from\kappa\to V_\kappa$ such
that for any $A\in H(\kappa^\plus)$ there is an embedding $j:V\to M$ having critical point $\kappa$ with
$j(\ell)(\kappa)=A$. This and the similar principles formulated for many other large cardinal notions,
including weakly compact, indescribable, unfoldable, Ramsey, strongly unfoldable and strongly compact
cardinals, are best conceived as an expression of the Laver function concept from supercompact cardinals
for these weaker large cardinal notions. The resulting Laver diamond principles $\LDstar_\kappa$ can hold
or fail in a variety of interesting ways.

The classical diamond principle $\Diamond_\kappa$, for an infinite cardinal $\kappa$, is a gem of modern
infinite combinatorics; its reflections have illuminated the path of innumerable combinatorial
constructions and unify an infinite amount of combinatorial information into a single, transparent
statement. When $\kappa$ exhibits any of a variety of large cardinal properties, however, the principle can
be easily strengthened in natural ways, and in this article I introduce and survey the resulting class of
strengthened principles, which I call the Laver diamond principles.

The classical diamond principle $\Diamond_\kappa$ asserts that there is a function $d\from\kappa\to
V_\kappa$ such that for any set $A\of\kappa$ the set $B=\set{\alpha\in\dom(d)\st
d(\alpha)=A\intersect\alpha}$ is stationary. Succinctly, the function $d$ anticipates every set
$A\of\kappa$ on a stationary set. The strong diamond principles I have in mind effectively ask for much
more than the stationarity of $B$. For example, if $\kappa$ is a measurable cardinal, then we might
naturally ask that the set $B$ has normal measure one. If it does, and $j:V\to M$ is the corresponding
ultrapower embedding by this normal measure, then a simple ultrapower calculation verifies that
$j(d)(\kappa)=A$. Thus, the function $d$ becomes a kind of Laver function for measurability, a function for
which the value of $j(d)(\kappa)$ can be arbitrarily prescribed. And since (as I will show) a measurable
cardinal can but need not have such a Laver function, I view the existence of such a function as a new
large cardinal combinatorial principle generalizing $\Diamond_\kappa$. Specifically, I define the Laver
diamond principle $\LDmeas_\kappa$ to assert that there is a function $\ell\from\kappa\to V_\kappa$ such
that every set $A\in H(\kappa^\plus)$ has an embedding $j:V\to M$ with $j(\ell)(\kappa)=A$. Similar
principles asserting the existence of an appropriate Laver function can be formulated for many other large
cardinal contexts.

The Laver function concept was first introduced by Richard Laver in his celebrated paper \cite{Laver78},
where he proved that every supercompact cardinal $\kappa$ has a function $\ell\from\kappa\to V_\kappa$ (now
known as a {\df Laver function}) such that for any set $A$ and any $\theta\geq|\TC(A)|$, there is a
$\theta$-supercompactness embedding $j:V\to M$ such that $j(\ell)(\kappa)=A$. Thus, in the supercompactness
context, a Laver function can anticipate any set at all, with an embedding as high in the supercompactness
hierarchy as desired. Various authors have generalized the Laver function concept to other large cardinal
contexts, for which corresponding functions can be proved to exist. Gitik and Shelah \cite{GitikShelah89},
for example, proved that every strong cardinal has a kind of Laver function, and Corazza \cite{CorazzaJSL}
showed the same for extendible cardinals (in the context of a general Laver function analysis). Since these
accounts have required that a Laver function anticipate every set $A$ in the universe, thereby implying
that $\kappa$ is a strong cardinal, the analysis has been mostly contained in the upper part of the large
cardinal hierarchy. In \cite{Hamkins2000:LotteryPreparation} I relaxed this requirement to a class of sets
$A$ appropriate for strongly compact embeddings (namely, the sets $A$ which appear in strongly compact
ultrapowers), and was able to force the existence of a strongly compact Laver function for any strongly
compact cardinal. In that article, I developed a general framework for analyzing generic Laver functions.
Lower down in the large cardinal hierarchy, Shelah and V\"a\"an\"anen
\cite{ShelahVaananan726:ExtensionsOfInfinitaryLogic} have considered independently a statement that is
equivalent to the Laver diamond for weakly compact cardinals. So the time may be ripe to consider the full
spectrum of the Laver function concept.

The possibility that there might not be Laver functions for certain large cardinals leads one naturally to
view the existence of a Laver function as a combinatorial principle of its own. This principle has been
particularly interesting for the smaller large cardinals, such as weakly compact, indescribable, unfoldable
or measurable cardinals, where it can hold or fail in a variety of ways.

Throughout this article, I use the three-dot notation $f\from A\to B$ to mean that $f$ is a partial
function from $A$ to $B$, that is, that $f$ is a function for which $\dom(f)\of A$ and $\ran(f)\of B$. I
denote by $f\image X$ the image of $X$ under $f$, that is, the set $\set{f(a)\st a\in X}$. In several
arguments below, I recursively build a partial function $\ell\from\kappa\to V_\kappa$ by first supposing
$\ell\restrict\gamma$ is defined, and then defining $\ell(\gamma)$ or throwing $\gamma$ out of
$\dom(\ell)$. Since the final function $\ell$ is only a partial function, this formalization is somewhat
loose, but I hope the reader can see that I am actually recursively building the (total) map $\gamma\mapsto
\ell\restrict\gamma$, and thereby avoid confusion. I use the term {\df $\leqkappa$-distributive} forcing to
refer to forcing notions that add no new sequences over the ground model of length less than or equal to
$\kappa$.

\Section The Laver diamond principles $\LaverDiamond_\kappa$

Let me now define the Laver diamond principles. The basic idea is that for each large cardinal notion
$\star$ and cardinal $\kappa$ exhibiting it, the principle $\LDstar_\kappa$ asserts that there is a
function $\ell\from\kappa\to V_\kappa$ such that for every (appropriate) set $A$, there is an (appropriate)
embedding $j$ such that $j(\ell)(\kappa)=A$. What is appropriate is determined by the large cardinal notion
under consideration. The $\LD_\kappa$ symbol was chosen for its vague resemblances both to
$\Diamond_\kappa$ and to the letter L, and it should be pronounced as ``{\it the Laver diamond principle}''
at $\kappa$.

I have mentioned already the principle $\LDmeas_\kappa$ in the case of a measurable cardinal $\kappa$,
which is a very natural case. A cardinal $\kappa$ is {\df measurable} if and only if there is a transitive
class $M$ and an elementary embedding $j:V\to M$ with critical point $\kappa$. For any such embedding, the
models $M$ and $V$ must agree on $H(\kappa^\plus)$, but not necessarily on larger sets. This suggests an
appropriate class of sets and embeddings, and the Laver diamond principle $\LDmeas_\kappa$ accordingly
asserts that there is a partial function $\ell\from\kappa\to V_\kappa$ such that for any set $A\in
H(\kappa^\plus)$ there is an elementary embedding $j:V\to M$ with critical point $\kappa$ such that
$j(\ell)(\kappa)=A$. If there is such an embedding, then it is easy to see that the ultrapower by the
induced normal measure also has this feature, and so one can assume without loss of generality that $j$ is
the ultrapower by a normal measure.

Let me consider the remaining large cardinal notions roughly in the order of their consistency strength,
beginning with the weakly compact cardinals. A cardinal $\kappa$ is {\df weakly compact} if for every
transitive structure $M$ of size $\kappa$ with $\kappa\in M$ there is an elementary embedding $j:M\to N$
into another transitive structure $N$ with critical point $\kappa$. (It suffices to consider only the
structures $M$ that are {\df nice} in the sense that $M$ satisfies some fragment of set theory and
$M^\ltkappa\of M$ as well.) The Laver diamond principle $\LDwc_\kappa$ for weak compactness is the
assertion that there is a function $\ell\from\kappa\to V_\kappa$ such that for any set $A\in
H(\kappa^\plus)$ and any transitive structure $M$ of size $\kappa$ with $A,\ell\in M$, there is a
transitive set $N$ and an elementary embedding $j:M\to N$ with critical point $\kappa$ such that
$j(\ell)(\kappa)=A$.

By \cite{Hauser1991:IndescribableCardinals}, a cardinal $\kappa$ is {\df $\Pi^m_n$-indescribable} if for
any nice structure $M$ of size $\kappa$ there is a transitive set $N$ and an elementary embedding $j:M\to
N$ with critical point $\kappa$ such that $N$ is $\Sigma^m_n$-correct, that is, $(V_{\kappa+m})^N\elesub_n
V_{\kappa+m}$ and $M^\card{V_{\kappa+m-2}}\of M$ (meaning $M^\ltkappa\of M$ when $m=1$). The principle
$\LaverDiamond_\kappa^{\hbox{\!\!\tiny$\Pi^m_n$-ind}}$ is the assertion that there is a function
$\ell\from\kappa\to V_\kappa$ such that for any set $A\in V_{\kappa+m}$ and any nice structure $M$ with
$\ell\in M$, there is a transitive set $N$ and an elementary embedding $j:M\to N$ with critical point
$\kappa$ such that $N$ is $\Sigma^m_n$-correct and $j(\ell)(\kappa)=A$.

A cardinal $\kappa$ is {\df unfoldable} if it is $\theta$-unfoldable for every ordinal $\theta$, that is,
if for every nice structure $M$ of size $\kappa$ there is a transitive structure $N$ and an elementary
embedding $j:M\to N$ such that $j(\kappa)>\theta$. The principle $\LDunf_\kappa$ is the assertion that
there is a function $\ell\from\kappa\to V_\kappa$ such that for any ordinal $\theta$, any set $A\in
H(\kappa^\plus)$ and any transitive set $M$ of size $\kappa$ with $M^\ltkappa\of M$ and $\ell,A\in M$,
there is a $\theta$-unfoldability embedding $j:M\to N$ with $j(\ell)(\kappa)=A$. The principle
$\LDthetaunf_\kappa$ is the assertion that there is a Laver function that works for $\theta$-unfoldability
embeddings. A function $\ell\from\kappa\to V_\kappa$ is an {\df ordinal anticipating} Laver function for
$\theta$-unfoldability if for any $\alpha<\theta$ and whenever $M$ is a nice structure containing $\ell$,
there is a $\theta$-unfoldability embedding $j:M\to N$ such that $j(\ell)(\kappa)=\alpha$ (that is, we
require $\ell$ to anticipate ordinals, but not necessarily other sets). Combining the two ideas, it would
be natural to inquire about a function $\ell$ that anticipates every set in $L_\theta(P(\kappa))$.

A cardinal $\kappa$ is {\df strongly unfoldable} if it is strongly $\theta$-unfoldable for every ordinal
$\theta$, that is, if for every nice structure $M$ there is a transitive structure $N$ and an elementary
embedding $j:M\to N$ with $V_\theta\of N$. The principle $\LDsunf_\kappa$ is the assertion that there is a
function $\ell\from\kappa\to V_\kappa$ such that for any cardinal $\theta$, any set $A\in H(\theta^\plus)$
and any transitive set $M$ of size $\kappa$ with $M^\ltkappa\of M$ and $\ell\in M$, there is a strong
$\theta$-unfoldability embedding $j:M\to N$ with $j(\ell)(\kappa)=A$. One obtains the principle
$\LDthetasunf_\kappa$ by fixing $\theta$.

For Ramsey cardinals, I use an embedding characterization mentioned to me by Philip Welch
\cite{Welch2002:PersonalCommunication}. Specifically, a cardinal $\kappa$ is a {\df Ramsey} cardinal if and
only if for any set $B\of\kappa$ there is a nice, iterable, weakly amenable structure $\<M,F>$ containing
$B$; that is, there is a structure $\<M,F>$ such that $M$ is a transitive set of size $\kappa$,
$M\satisfies\ZFC^-$, $M^\ltkappa\of M$, $\kappa,B\in M$, and $F$ is an iterable weakly amenable $M$-normal
measure, meaning that $F\intersect D\in M$ for any $D\in H(\kappa^\plus)^M$, $F$ is closed under diagonal
intersections of $\kappa$-sequences in $M$, and all iterations of the ultrapower of $M$ by $F$ are
well-founded. The principle $\LDram_\kappa$ is the assertion that there is a function $\ell\from\kappa\to
V_\kappa$ such that for any set $A\in H(\kappa^\plus)$ and any $B\of\kappa$ there is such a structure
$\<M,F>$ as above with $\ell,A\in M$ also, and the ultrapower of $M$ by $F$ has $j(\ell)(\kappa)=A$. Using
the notion of iterable extender embeddings, one could define a notion of {\df ordinal anticipating} Laver
functions for Ramsey cardinals, by insisting that the function $\ell$ anticipate all ordinals via such
embeddings.

A cardinal $\kappa$ is {\df strong} if it is $\theta$-strong for every ordinal $\theta$, that is, if there
is a transitive class $M$ and an embedding $j:V\to M$ with critical point $\kappa$ such that $V_\theta\of
M$. The standard analysis of such embeddings (e.g. \cite{Kanamori:TheHigherInfinite}) shows that if there
is such an embedding, then there is one that is that is the ultrapower by an {\df extender}, a certain
directed system of ultrafilters, and so the notion of strongness is a first-order concept in set theory.
The principle $\LDstr_\kappa$ is the assertion that there is a function $\ell\from\kappa\to V_\kappa$ such
that for any $\theta$ and any $A\in V_\theta$ there is a $\theta$-strongness embedding $j:V\to M$ such that
$j(\ell)(\kappa)=A$. The local principle $\LDthetastr_\kappa$ makes the same assertion when $\theta$ is
fixed.

A cardinal $\kappa$ is {\df superstrong} (with target $\lambda$) when it is the critical point of an
embedding $j:V\to M$ such that $j(\kappa)=\lambda$ and $V_\lambda\of M$. The principle
$\LD_\kappa^{\hbox{\!\!\tiny superstrong}}$ (with target $\lambda$) is the assertion that there is a
function $\ell\from\kappa\to V_\kappa$ such that for any $A\in V_\lambda$ there is a superstrong embedding
$j:V\to M$ with target $\lambda$ such that $j(\ell)(\kappa)=A$.

A cardinal $\kappa$ is {\df strongly compact} if it is $\theta$-strongly compact for every cardinal
$\theta$, that is, if there is a transitive class $M$ and an embedding $j:V\to M$ such that $M$ has the
{\df $\theta$-cover property}, that is, for every set $X\of M$ of size at most $\theta$ there is a set
$Y\in M$ with $X\of Y$ and $|Y|^M<j(\kappa)$. This property follows if $j$ is an ultrapower embedding by a
measure on some set and there is a set $s\of j(\kappa)$ in $M$ of size less than $j(\kappa)$ there with
$j\image\theta\of s$. It is easy to verify that the measure induced by such a set $s$ is a fine measure on
$P_\kappa\theta$ (in the nontrivial case that $\theta\geq\kappa$), and conversely, the ultrapower by any
fine measure on $P_\kappa\theta$ gives rise to such an embedding. In general, one can't expect agreement
between $M$ and $V$ beyond $H(\kappa^\plus)$, but $j(\kappa)$ must be above $\theta$. The principle
$\LDstrc_\kappa$ is the assertion that there is a function $\ell\from\kappa\to V_\kappa$ such that for any
cardinal $\theta$ and any set $A\in H(\kappa^\plus)$, there is a $\theta$-strong compactness embedding
$j:V\to M$ with $j(\ell)(\kappa)=A$. The local principle $\LDthetastr_\kappa$ makes this assertion when
$\theta$ is fixed. An {\df ordinal anticipating} Laver function must anticipate all ordinals
$\alpha<\theta$ (but not necessarily other sets).

A cardinal $\kappa$ is {\df supercompact} if it is $\theta$-supercompact for every cardinal $\theta$, that
is, if there is a transitive class $M$ and an embedding $j:V\to M$ with critical point $\kappa$ such that
$M^\theta\of M$ and $j(\kappa)>\theta$. In particular, $j\image\theta\in M$. It is easy to see that the
measure induced by $j\image\theta$ is a normal fine measure on $P_\kappa\theta$, and conversely, the
ultrapower by any normal fine measure on $P_\kappa\theta$ gives rise to a $\theta$-supercompactness
embedding. Because $M^\theta\of M$, it follows that $M$ and $V$ agree on $H(\theta^\plus)$. The principle
$\LDsc_\kappa$ is the assertion that there is a function $\ell\from\kappa\to V_\kappa$ such that for any
$\theta$ and any $A\in H(\theta^\plus)$ there is a $\theta$-supercompactness embedding $j:V\to M$ such that
$j(\ell)(\kappa)=A$. The local principle $\LDthetasc_\kappa$ makes this assertion for fixed $\theta$.

A cardinal $\kappa$ is {\df extendible} if for every ordinal $\alpha>\kappa$ there is an ordinal $\beta$
and an embedding $j:V_\alpha\to V_\beta$ with critical point $\kappa$. The principle $\LDext_\kappa$ is the
assertion that there is a function $\ell\from\kappa\to V_\kappa$ such that for every set $A$ and any
ordinal $\alpha>\kappa$ there is an embedding $j:V_\alpha\to V_\beta$ such that $j(\ell)(\kappa)=A$.

A cardinal $\kappa$ is {\df almost huge} with target $\lambda$ (or {\df almost $\lambda$-huge}) if there is
an embedding $j:V\to M$ with critical point $\kappa$ such that $\lambda=j(\kappa)$ and $M^\ltlambda\of M$.
The cardinal $\kappa$ is {\df super almost huge} if it is almost huge with arbitrarily large targets
$\lambda$. The principle $\LDlambdaahuge_\kappa$ is the assertion that there is a function
$\ell\from\kappa\to V_\kappa$ such that for any $A\in V_\lambda$ there is an almost hugeness embedding
$j:V\to M$ with target $\lambda$ such that $j(\ell)(\kappa)=A$. The principle $\LDsahuge_\kappa$ makes this
assertion about $\ell$ for arbitrarily large targets $\lambda$. By increasing the degree of closure to
$M^\lambda\of M$, one obtains the notion of $\kappa$ being {\df huge} with target $\lambda$, the notion of
$\kappa$ being {\df superhuge}, and the principles $\LDlambdahuge_\kappa$ and $\LDshuge_\kappa$. By
increasing the degree of closure to $M^{j^n(\kappa)}\of M$, where $n\in \omega$, one obtains the notions of
{\df $n$-huge} with target $\lambda$ and {\df super $n$-huge}, with the corresponding Laver diamond
principles $\LDlambdanhuge_\kappa$ and $\LDsnhuge_\kappa$.

I trust that the reader will be able to extend the definition of the Laver diamond principle to any large
cardinal context that has a natural class of embeddings $j$ and a natural class of sets $A$ to be
anticipated.

\Section Automatic instances of the Laver diamond principles $\LaverDiamond_\kappa$

Certain instances of the Laver diamond principle are provable in \ZFC. Let us begin with Laver's
\cite{Laver78} celebrated result.

\Theorem.(Laver) If\/ $\kappa$ is supercompact, then $\LDsc_\kappa$.\label{Laver}

\Proof: I sketch the proof of \cite{Laver78}, constructing the function $\ell$ by transfinite recursion.
If\/ $\ell\restrict\gamma$ has been defined, then let $\lambda$ be least such that there is some set $a\in
H(\lambda^\plus)$ which is not anticipated by any $\lambda$-supercompactness embedding $j:V\to M$ with
critical point $\gamma$, that is, such that $j(\ell\restrict\gamma)(\gamma)\not=a$ for all such $j$. If
there is such a $\lambda$, then let $\ell(\gamma)=a$ for some such $a$. A simple reflection argument shows
that $\lambda$ must be below the next strong cardinal if it exists at all. I claim that $\ell$ is a
supercompactness Laver function. If not, there is some least $\theta$ and some $a\in H(\theta^\plus)$ which
is not anticipated by $\ell$. Fix any $2^{\theta^\ltkappa}$-supercompactness embedding $j:V\to M$. Note
that $M$ and $V$ have the same fine normal measures on $P_\kappa\theta$ and the ultrapower embeddings by
these measures in $V$ or $M$ will agree on $\ell$. Thus, the model $M$ agrees that $\theta$ is least such
that some $a\in H(\theta^\plus)$ is not anticipated by $\ell=j(\ell)\restrict\kappa$, and so
$j(\ell)(\kappa)=a$ for some such $a$. Let now $j_0$
\begin{figure}[h]\factordiagramup{V}{j}{M}{j_0}{k}{M_0}\end{figure} be the $\theta$-supercompactness factor
embedding induced by $j$. Specifically, let $X=\set{j(f)(j\image\theta)\st f\in V}$ and observe that
$X\elesub M$ by the Tarski-Vaught criterion. Let $\pi:X\iso M_0$ be the Mostowski collapse of\/ $X$, and
since $\ran(j)\of X$, let $j_0=\pi\compose j$. This produces the commutative diagram, where $k$ is the
inverse of $\pi$. Applying $\pi$ to the definition of\/ $X$ yields that
$M_0=\set{j_0(f)(j_0\image\theta)\st f\in V}$. It follows that $j_0$ is the ultrapower of\/ $V$ by the
normal fine measure $\mu_0=\set{A\of P_\kappa\theta\st j_0\image\theta\in j_0(A)}$; the mapping
$j_0(f)(j_0\image\theta)\mapsto [f]_{\mu_0}$ is an isomorphism of\/ $M_0$ to $\Ult(V,\mu_0)$. Thus, the set
$a$ is anticipated by a normal fine measure on $P_\kappa\theta$ after all, a contradiction.\QED

I invite the reader to compare the proof just given to the usual proof of $\Diamond$ in $L$; at each stage
$\gamma$, one defines $\ell(\gamma)$ to be a least counterexample, if any exists. At the top, then, there
can be no counterexample at $\kappa$, because any value of $j(\ell)(\kappa)$ in a larger embedding will
betray itself in the induced factor embedding.

Because the definition of $\ell$ in the previous theorem was local---the value of $\ell(\gamma)$ depended
essentially only on $\ell\restrict\gamma$---the proof actually establishes a universal kind of Laver
function. Specifically, let me define $\LDsc$ to be assertion that there is a (possibly proper) class
function $\ell\from\ORD\to V$, called a {\df universal} Laver function, such that whenever $\kappa$ is a
supercompact cardinal, then $\ell\restrict\kappa$ witnesses $\LDsc_\kappa$. And there are analogous
universal forms of the other Laver diamond principles, such as $\LDwc$, $\LDunf$, $\LDmeas$, and so on (but
unlike $\LDsc$, these are not generally theorems of \ZFC). Such a universal Laver function weaves together
a single coherent sequence from the possible Laver functions at each of the various cardinals under
consideration.  Such a weaving together of the Laver functions seems non-trivial especially at compound
limits of the cardinals, such as supercompact limits of supercompact cardinals, since one must in a sense
hope that the function constructed up to that stage works at that stage. Nevertheless, since in the proof
of Theorem \ref{Laver} the definition of $\ell(\gamma)$ did not depend on $\kappa$, the final length of the
Laver function, we do indeed obtain a universal function (and this is observed also in
\cite{KimchiMagidor}):

\Theorem. (Assuming AC for classes) The principle $\LDsc$ holds. That is, there is a universal Laver
function $\ell\from\ORD\to V$ such that $\ell\restrict\kappa$ is a $\LDsc_\kappa$ Laver function for any
supercompact cardinal $\kappa$.\label{UniformLDsc}

One uses the principle of AC for classes to obtain a well-ordering $\trianglelt$ of $V$ in the case that
there are a proper class of supercompact cardinals. This hypothesis is harmless; it is conservative over
\ZFC\ for purely first order statements, since one can force it over any model of set theory without adding
sets. Furthermore, if there are a proper class of supercompact cardinals, then this assumption is
necessary, because the range of $\ell$ will necessarily include all sets, and from $\ell$ one can define a
well-ordering of $V$.

One cannot expect to obtain a universal Laver function by simply pasting together Laver functions for
smaller cardinals, in view of the following.

\Observation. If $\kappa$ is the least supercompact limit of supercompact cardinals, then there is a
function $\ell\from\kappa\to V_\kappa$ such that $\ell\restrict\delta$ is a $\LDsc_\delta$ Laver function
for every supercompact cardinal $\delta<\kappa$, but $\ell$ is not a $\LDsc_\kappa$ Laver function.

\Proof: Let $\bar\ell$ be a $\LDsc_\kappa$ Laver function for $\kappa$, as constructed in Theorems
\ref{Laver} or \ref{UniformLDsc}, and define $\ell\from\kappa\to V_\kappa$ by
$\ell(\gamma)=\bar\ell(\gamma)$ unless $\gamma$ is a limit of supercompact cardinals, in which case let
$\ell(\gamma)=17$. I claim that if $\delta<\kappa$ is a supercompact cardinal, then $\ell\restrict\delta$
is a $\LDsc_\delta$ Laver function, because $\delta$ is not a limit of supercompact cardinals ($\kappa$
being the least such supercompact cardinal), and so the change from $\bar\ell$ does not affect
$j(\bar\ell)(\delta)=j(\ell)(\delta)$ for any embedding $j:V\to M$ with critical point $\delta$. That is,
$\bar\ell$ and $\ell$ agree on all the points that are relevant to any supercompact cardinal
$\delta<\kappa$, because such a $\delta$ will never be a limit of supercompact cardinals. But at the top,
if $j:V\to M$ is any embedding having critical point $\kappa$, then necessarily $j(\ell)(\kappa)=17$, and
so $\ell$ is definitely not a $\LDsc_\kappa$ function.\QED

Similarly, one should not expect to prove in general that if $\ell$ is a $\LDsc_\kappa$ Laver function for
a supercompact cardinal $\kappa$, then $\ell\restrict\gamma$ is a $\LDmeas_\gamma$ Laver function for every
measurable cardinal $\gamma<\kappa$, because in witnessing $\LDsc_\kappa$ one only needs to consider
supercompactness embeddings $j:V\to M$ of very high Mitchell rank, and so in particular they all have the
feature that $\kappa$ is a measurable cardinal of very high Mitchell rank in $M$. The function $\ell$,
therefore, needs to be defined only on such cardinals. In particular, if we trivialize the values of $\ell$
elsewhere, then $\ell$ will still witness $\LDsc_\kappa$, but if $\gamma<\kappa$ is measurable with low
Mitchell rank, then for any embedding $j:V\to M$ with critical point $\gamma$, the value of
$j(\ell)(\gamma)$ will be trivial. So $\ell\restrict\gamma$ will not be a $\LDmeas_\gamma$ Laver function.
Similar observations apply to the partially supercompact cardinals below $\kappa$.

One can, of course, show that if $\ell$ is a $\LDsc_\kappa$ Laver function (or much less, merely a
$\LD^{\hbox{\!\!\tiny$o(\kappa){\geq}1$}}_\kappa$ Laver function), then nevertheless $\ell\restrict\gamma$
will witness $\LDmeas_\gamma$ for many $\gamma<\kappa$. An easy reflection argument shows that this is true
on a set of normal measure one. Conversely, if $\gamma$ is a measurable cardinal with nontrivial Mitchell
rank and $\LDmeas_\gamma$ holds for all $\gamma<\kappa$, then $\LDmeas_\kappa$ also holds, since it holds
in $M$ for some embedding $j:V\to M$ in which $\kappa$ is measurable, and it is absolute from $M$ to $V$.

The question of whether restrictions of a given Laver function witness lesser Laver diamond principles
below arises again below in Theorem \ref{SimultaneousLD}, where I show how to produce universal Laver
functions that simultaneously witness the Laver diamond for several different large cardinal notions.

The proof of Theorem \ref{Laver} can be refined into the supercompactness hierarchy. In order to get a
Laver function for $\theta$-supercompactness embeddings, we didn't need to know that $\kappa$ was fully
supercompact, but rather only that it was $2^{\theta^\ltkappa}$-supercompact. Thus, Laver's argument
establishes:

\Corollary. If\/ $\kappa$ is $2^{\theta^\ltkappa}$-supercompact, then $\LDthetasc_\kappa$ holds.
\label{LaverLocal}

And this reasoning shows that the construction of a universal Laver function $\ell\from\ORD\to V$ for
$\LDsc$ in Theorem \ref{UniformLDsc} actually serves as a Laver function for many partially supercompact
cardinals. Specifically, if $\kappa$ is $2^{\theta^\ltkappa}$-supercompact, then the argument shows that
$\ell\restrict\kappa$ (or more properly, $\ell\intersect\kappa\cross V_\kappa$, because the range should be
in $V_\kappa$) will be a $\LDthetasc_\kappa$ Laver function.

It is natural to inquire whether the hypothesis of Corollary \ref{LaverLocal} can be reduced to the
assumption only that $\kappa$ is $\theta$-supercompact, in an attempt to match the large cardinal strength
of the conclusion with that of the hypothesis.

\Question. If $\kappa$ is $\theta$-supercompact and $\kappa<\theta$, must $\LDthetasc_\kappa$
hold?\label{LDthetascQuestion}

I conjecture a negative answer (meaning that I expect that a negative answer is relatively consistent). In
the case $\theta=\kappa$ this is already known, because $\kappa$ is $\kappa$-supercompact if and only if
$\kappa$ is measurable and $\LDkappasc_\kappa$ is the same as $\LDmeas_\kappa$, so Theorem
\ref{LDmeasFailsInL[mu]} shows that $\kappa$ can be $\kappa$-supercompact even when $\LDkappasc_\kappa$
fails. But the question is open at higher levels of supercompactness, largely because (1) we have no inner
model theory for such cardinals, where one might expect $\LDthetasc_\kappa$ to fail, and (2) we have not
been able to force $\neg\LDthetasc_\kappa$ while preserving the $\theta$-supercompactness of $\kappa$.
Despite my conjecture that we will eventually build a model of $\neg\LDthetasc_\kappa$ (probably by
forcing), Corollary \ref{ForcingLDscLocal} does show that in terms of consistency strength, the hypothesis
of Corollary \ref{LaverLocal} can be reduced, for whenever $\kappa$ is $\theta$-supercompact, then
$\LDthetasc_\kappa$ is forceable. The argument of Theorem \ref{FewMeasures} shows the connection between
Question \ref{LDthetascQuestion} and the question of the number of normal fine measures on $P_\kappa\theta$
that a $\theta$-supercompact cardinal $\kappa$ can or must have.

Perhaps Corollary \ref{LaverLocal} reveals the reason why Laver functions can be proved to exist for
supercompact cardinals: it is simply because supercompactness comes in a hierarchy for which the existence
of Laver functions at each level is implied by the higher levels. The same phenomenon will occur in Theorem
\ref{GitikShelah} below with strong cardinals, and in Theorem \ref{Corazza} with superhuge cardinals. But
it does not occur with weakly compact, indescribable, unfoldable and measurable cardinals.

Let me turn now to the case of strong cardinals.

\Theorem.(Gitik-Shelah \cite{GitikShelah89}) If\/ $\kappa$ is a strong cardinal, then
$\LDstr_\kappa$.\label{GitikShelah}

\Proof: This is similar to Theorem \ref{Laver}. Suppose we have defined $\ell$ up to $\gamma$, and we want
to define $\ell(\gamma)$. Let $\lambda$ be least such that there is some $a\in V_\lambda$ such that no
$\lambda$-strongness embedding $j:V\to M$ with critical point $\gamma$ has $j(\ell)(\gamma)=a$. Define
$\ell(\gamma)=a$ for some such $a$. We claim that $\ell$ is a $\LDstr_\kappa$ Laver function. If not, then
let $\theta$ be least such that some $a\in V_\theta$ is not anticipated by $\ell$. Fix any
$(\theta+1)$-strongness embedding $j:V\to M$ with critical point $\kappa$. Since $M$ and $V$ have the same
$V_{\theta+1}$, they have the same extenders for $\theta$-strongness embeddings, and so $M$ agrees that
$\theta$ is least such that there is some $a$ that is not anticipated by $\ell=j(\ell)\restrict\kappa$.
Thus, $j(\ell)(\kappa)=a$ for some such $a$. Let $j_0$ be the induced factor embedding induced by the
$\theta$-strongness extender. \factordiagramup{V}{j}{M}{j_0}{k}{M_0} Since the critical point of $k$ is at
least $\theta$, it follows as earlier that $j_0(\ell)(\kappa)=a$. Thus, the set $a$ is anticipated by a
$\theta$-strongness embedding after all, a contradiction.\QED

This theorem also admits universal and local forms, as we obtained above in the case of supercompact
cardinals.

\Corollary. \
\begin{enumerate}\item (Assuming AC for classes) There is a universal $\LDstr$ Laver function.
\item If\/ $\kappa$ is $(\theta+1)$-strong, then $\LDthetastr_\kappa$ holds. Indeed, the restriction to
$\kappa$ of the $\LDstr$ Laver function constructed in Theorem \ref{GitikShelah} witnesses
$\LDthetastr_\kappa$.
\end{enumerate}\label{UniformLDstr}

In many cases one can combine the various constructions to produce a single universal Laver function
$\ell\from\ORD\to V$ that simultaneously witnesses the Laver Diamond for many different large cardinal
notions. Let me give the first hint of this in the following theorem.

\Theorem. (Assuming AC for classes) There is a function $\ell\from\ORD\to V$ which simultaneously witnesses
$\LDsc$ and $\LDstr$. In addition, if $\kappa$ is strong and $2^{\theta^\ltkappa}$-supercompact, where
$\kappa\leq\theta$, then $\ell\restrict\kappa$ is a $\LDthetasc_\kappa$ Laver function. And if $\kappa$ is
$\theta+1$-strong, then $\ell\intersect\kappa\cross V_\kappa$ is a $\LDthetastr_\kappa$ Laver
function.\label{SimultaneousLD}

\Proof: The point is that the domains of the $\LDsc$ and $\LDstr$ Laver functions constructed in Theorems
\ref{Laver} and \ref{GitikShelah} need not interfere with each other, and so one can simply take the union
of the two functions.

\SubLemma. The restriction of any $\LDsc$ Laver function to the class of non-supercompact strong cardinals
still witnesses $\LDsc$.\label{LDscDomain}

\Proof: Suppose that $\ell$ witnesses $\LDsc$, and let $\ell^*$ be the restriction of $\ell$ to the
non-supercompact strong cardinals. I claim by induction that $\ell^*\restrict\kappa$ witnesses
$\LDsc_\kappa$ for any supercompact cardinal. If $\kappa$ is supercompact, then for any $\theta\geq
2^\kappa$ and any $A\in H(\theta^\plus)$ there is a $\theta$-supercompactness embedding $j:V\to M$ such
that $j(\ell)(\kappa)=A$. First, I observe that since $\theta\geq 2^\kappa$, the induced extender
$j\restrict P(\kappa)$ is in $M$, and consequently $\kappa$ is superstrong with target $j(\kappa)$ in $M$.
Since $j(\kappa)$ is supercompact in $M$, it follows that $\kappa$ is a strong cardinal in $M$. Next, if
$\kappa$ is not supercompact in $M$, then $j(\ell^*)(\kappa)=A$, and we are done. Otherwise, $\kappa$ is
supercompact in $M$, and by the induction hypothesis $j(\ell^*)\restrict\kappa=\ell^*\restrict\kappa$
witnesses $\LDsc_\kappa$ in $M$. Up to $\theta$-supercompactness, this is absolute to $V$, since any
$\theta$-supercompactness measure in $M$ is a $\theta$-supercompactness measure in $V$, and the operation
of such measures on $P(\kappa)$ in $M$ is the same as in $V$. Thus, there is a $\theta$-supercompactness
embedding $j_0:V\to M_0$ such that $j(\ell^*)(\kappa)=A$, as desired.\QED

\SubLemma. The restriction of any $\LDstr$ Laver function to the class of non-strong measurable cardinals
still witnesses $\LDstr$.\label{LDstrDomain}

\Proof: I argue similarly here. Suppose that $\ell$ witnesses $\LDstr$, and let $\ell^*$ be the restriction
of $\ell$ to the non-strong measurable cardinals. If $\kappa$ is strong and $A\in V_\theta$, where
$\kappa+2\leq\theta$, then there is a $\theta$-strongness embedding $j:V\to M$ such that
$j(\ell)(\kappa)=A$. Since $\kappa+2\leq\theta$, it follows that $\kappa$ is measurable in $M$. If $\kappa$
is not strong in $M$, then $j(\ell^*)(\kappa)=A$, and we are done. Otherwise, $\kappa$ is strong in $M$,
and so by induction we know that $j(\ell^*)\restrict\kappa=\ell^*\restrict\kappa$ witnesses $\LDstr_\kappa$
in $M$. So there is a $\theta$-strongness extender in $M$ with $\ell^*$ anticipating $A$. Since
$V_\theta\of M$, this extender may be applied in $V$ to obtain $\theta$-strongness embedding $j_0:V\to M_0$
with $j_0(\ell^*)(\kappa)=A$, as desired.\QED

Combining these two facts, one simply takes the union of the two functions constructed in Theorem
\ref{UniformLDsc} and Corollary \ref{UniformLDstr}, restricted as in the lemmas, and this produces a
function $\ell$ witnessing both $\LDstr$ and $\LDsc$. Since every supercompact cardinal is a limit of
strong cardinals, it follows that every supercompact cardinal is closed under any $\LDstr$ Laver function,
in the sense that $\ell\image\kappa\of V_\kappa$. Similarly, if $\ell$ is the $\LDsc$ Laver function
constructed in Theorem \ref{UniformLDsc}, then we already mentioned that $\ell(\gamma)$ is below the next
strong cardinal; so every strong cardinal is closed under the resulting $\LDsc$ Laver function. Thus, every
supercompact or strong cardinal is closed under the union of the two functions.

Let me now prove the more local facts in the latter part of the theorem. If $\kappa$ is strong and
$2^{\theta^\ltkappa}$-supercompact, then for any witnessing embedding $j:V\to M$, it follows as above that
$\kappa$ is superstrong in $M$ with target $j(\kappa)$, which is strong in $M$, so $\kappa$ is strong in
$M$. Thus, with these embeddings, we are evaluating the $\LDsc$ Laver function at $\kappa$ in $M$, rather
than the $\LDstr$ function, which is not defined on strong cardinals, and so the union function will
witness $\LDthetasc_\kappa$, just as in Corollary \ref{LaverLocal}. Finally, if $\kappa$ is
$\theta+1$-strong, then $\ell$ might not be closed under $\kappa$, since some smaller $\ell(\gamma)$ could
jump over $\kappa$, as a result of the supercompactness Laver function. Nevertheless, the function
$\ell\intersect \kappa\cross V_\kappa$ will witness $\LDthetastr_\kappa$ just as in Corollary
\ref{UniformLDstr}.\QED

The two lemmas can easily be strengthened to limit the domain even further: for example, any $\LDsc$ Laver
function can be restricted to the class of non-supercompact strong cardinals $\kappa$ that are at least
$2^{2^\kappa}$-supercompact, and any $\LDstr$ Laver function can be restricted to the class of non-strong
cardinals $\kappa$ that are at least $(2^{\kappa^\plus}{+}17)$-strong. These strengthenings, however, have
consequences for witnessing the Laver diamond for weaker amounts of partial supercompactness (the last part
of Theorem \ref{SimultaneousLD}), so one needs to trade off these two notions, or else add in another Laver
function that works with the weaker cardinals.

The next theorem shows that this simultaneous Laver function idea can be pushed quite a long way.

\Theorem. If the principles $\LDmeas$ and $\LDwc$ hold, then there is a single class function
$\ell\from\ORD\to V$ simultaneously witnessing $\LDsc$, $\LDstr$, $\LDmeas$ and $\LDwc$.

\Proof:  The idea is that the two lemmas below show that the domains of the $\LDmeas$ and $\LDwc$ Laver
functions need interfere neither with each other nor with the $\LDsc$ and $\LDstr$ Laver functions.

\SubLemma. The restriction of any $\LDmeas$ Laver function to the class of non-measurable weakly compact
cardinals still witnesses $\LDmeas$.\label{LDmeasDomain}

\Proof: Suppose that $\ell$ is a universal $\LDmeas$ Laver function, and let $\ell^*$ be the restriction of
$\ell$ to the class of non-measurable weakly compact cardinals. I will argue that $\ell^*$ still witnesses
$\LDmeas$. Suppose, accordingly, that $\kappa$ is a measurable cardinal and $A\in H(\kappa^\plus)$, and
that $\ell^*\restrict\delta$ witnesses $\LDmeas_\delta$ for any measurable cardinal $\delta<\kappa$. Since
$\ell$ witnesses $\LDmeas$, there is an embedding $j:V\to M$ with critical point $\kappa$ and
$j(\ell)(\kappa)=A$. Certainly $\kappa$ is weakly compact in $M$. If $\kappa$ is not measurable in $M$,
then also $j(\ell^*)(\kappa)=A$, and we are done. Otherwise, $\kappa$ is measurable in $M$, and since
$\kappa<j(\kappa)$, we know by (applying $j$ to) the induction hypothesis that $j(\ell^*)\restrict\kappa$
witnesses $\LDmeas_\kappa$ in $M$. But $j(\ell^*)\restrict\kappa=\ell^*\restrict\kappa$, and the witness of
$\LDmeas_\kappa$ is absolute from $M$ to $V$, since they share the same $P(\kappa)$. So in any case, we're
done.\QED

\SubLemma. The restriction of any $\LDwc$ Laver function to the class of non-weakly compact cardinals still
witnesses $\LDmeas$.\label{LDwcDomain}

\Proof: Suppose that $\ell$ witnesses $\LDwc$, and let $\ell^*$ be the restriction of $\ell$ to the class
of non-weakly compact cardinals. If $\kappa$ is any weakly compact cardinal, I claim that
$\ell^*\restrict\kappa$ witnesses $\LDwc_\kappa$. To see this, fix any nice structure $M$ and $A\in
H(\kappa^\plus)$. Let $\bar M$ be a larger nice structure, such that $M\in \bar M$ and $\bar
M\satisfies|M|=\kappa$. Since $\ell$ witnesses $\LDwc$, there is $j:\bar M\to \bar N$ such that
$j(\ell\restrict\kappa)(\kappa)=A$. If $\kappa$ is not weakly compact in $\bar N$, then
$j(\ell^*\restrict\kappa)(\kappa)=A$ also, and by restricting to the embedding $j\restrict M:M\to j(M)$, we
are done. Otherwise, $\kappa$ is weakly compact in $\bar N$, and so by induction we know that
$j(\ell^*)\restrict\kappa=\ell^*\restrict\kappa$ witnesses $\LDwc_\kappa$ in $\bar N$. In particular,
inside $\bar N$ there is an embedding $j_0:M\to N_0$ with $j(\ell^*\restrict\kappa)(\kappa)=A$, and the
lemma is proved.\QED

One now simply takes the union of the corresponding Laver functions, and obtains a single Laver function
simultaneously witnessing the corresponding Laver diamond principles. So the theorem is proved.\QED

One can incorporate partially supercompact and partially strong cardinals into the previous argument, and
the lemmas generalize to unfoldable, strongly unfoldable and Ramsey cardinals, and so on. The point is that
by paying careful attention to the domain of the Laver functions, one can effectively combine them into a
single function that simultaneously witnesses many Laver diamond principles at once. The goal is to
construct a single function $\ell\from\ORD\to V$ that simultaneously witnesses as many of the diamond
principles as possible.

Like Lemmas \ref{LDscDomain} and \ref{LDstrDomain}, the results of Lemmas \ref{LDmeasDomain} and
\ref{LDwcDomain} apply only to the universal forms of the Laver diamond principles $\LDmeas$ and $\LDwc$,
rather than to assertions of the Laver diamond at a particular cardinal. In particular, one cannot hope to
prove that every $\LDwc_\kappa$ Laver function can be restricted to the non-weakly compact cardinals and
still witness $\LDwc_\kappa$, because by Observation \ref{TrivialObservation} and the subsequent remarks,
any $\LDmeas_\kappa$ Laver function also witness $\LDwc$, and such functions concentrate on the weakly
compact cardinals. Similarly, one cannot expect to prove that every $\LDmeas_\kappa$ Laver function can be
restricted to the non-measurable cardinals, because in Theorem \ref{SelfReflectiveMeasurable} below, the
domain of the $\LDmeas_\kappa$ Laver function can be taken to consist entirely of measurable cardinals.

\Question. If\/ $\LDmeas_\kappa$ holds, is there a $\LDmeas_\kappa$ Laver function defined only on
non-measurable cardinals? If\/ $\LDwc_\kappa$ holds, is there a $\LDwc_\kappa$ Laver function defined only
on non-weakly compact cardinals?

If the notation $\LDmeas_\kappa(D)$ means that there is a $\LDmeas_\kappa$ Laver function with domain $D$,
then the question asks whether $\LDmeas_\kappa$ is equivalent to $\LDmeas_\kappa(N)$ and whether $\LDwc$ is
equivalent to $\LDwc(M)$, where $N$ and $M$ are the set of non-measurable cardinals and the set of
non-weakly compact cardinals below $\kappa$, respectively. One has versions of this question for a great
many large cardinals.

The reflective construction of Theorems \ref{Laver} and \ref{GitikShelah} applies to a great many other
large cardinal notions. For example, Paul Corazza proved the following:

\Theorem.(Corazza \cite[Prop 4.5, Thm 5.3, 5.4]{CorazzaJSL}) \\ 1. If\/ $\kappa$ is extendible, then
$\LDext_\kappa$.\\ 2. If\/ $\kappa$ is slightly more than super almost huge,
   then $\LDsahuge_\kappa$.
\label{Corazza}

Let me now turn to much smaller cardinals. It turns out that one can also carry out the reflective
construction down low.

\Theorem. If\/ $\kappa$ is a measurable cardinal, then $\LDwc_\kappa$ holds.

\Proof: Suppose that $\kappa$ is a measurable cardinal, and that $\ell$ has been defined up to stage
$\gamma$, with $\ell\image\gamma\of V_\gamma$. If there is some transitive set $M$ of size $\gamma$ and a
set $A\in H(\gamma^\plus)^M$ with $\ell,\gamma\in M$ such that no weak compactness embedding $j:M\to N$ has
$j(\ell)(\gamma)=A$, then let $\ell(\gamma)=A$ for some such $A$; otherwise, if there is no such
problematic pair, then let $\ell(\gamma)$ be anything in $V_\kappa$ that you fancy. I claim that $\ell$
witnesses $\LDwc_\kappa$. If not, there is some $M$ of size $\kappa$ and set $A\in H(\kappa^\plus)^M$ with
$\ell,\kappa,A\in M$ but no embedding $j:M\to N$ with $j(\ell)(\kappa)=A$. Fix any ultrapower embedding
$j:V\to \bar M$ by a normal measure on $\kappa$. Since $\bar M$ and $V$ have the same $H(\kappa^\plus)$,
they agree that there is such a problematic pair $\<M,A>$. Thus, $j(\ell)(\kappa)=A$ for some such
problematic pair in $\bar M$. Now consider the mapping $j\restrict M:M\to j(M)$. This object has hereditary
size $\kappa$, and so it is in $\bar M$. Also, it clearly has $j(\ell)(\kappa)=A$. Thus, it shows that the
pair $\<M,A>$ is not a problem after all, a contradiction.\QED

The previous argument can be strengthened to the following.

\Theorem. If\/ $\kappa$ is a measurable cardinal, then $\LDunf_\kappa$ holds.\label{MeasurableLDunf}

\Proof: The construction will produce an ordinal anticipating Laver function. If $\ell$ has been defined up
to stage $\gamma$, let $\lambda$ be least, if any, such that for some $a=\alpha<\lambda$ or $a=A\in
H(\gamma^\plus)$ and some transitive set $N$ of size $\gamma$ with $\ell\restrict\gamma\in N$, there is no
$\lambda$-unfoldability embedding $j:N\to M$ having $j(\ell)(\gamma)=a$. Let $\ell(\gamma)=a$ for some such
$N$ and $a$. Necessarily $\lambda<\kappa$, if it exists, since if every $\alpha<\kappa$ is possible, then
for any $\theta$ fix any embedding $\bar j:V\to \bar M$ with critical point $\kappa$ and $\bar
j(\kappa)>\theta$ (simply iterate a normal measure on $\kappa$), and conclude that in $\bar M$ any
$\alpha<\theta$ and any $A\in H(\gamma^\plus)$ is possible with $N$ (contradicting the assumption that a
bad $\lambda$ exists). This defines $\ell\from\kappa\to V_\kappa$.

Suppose towards a contradiction that this is not an unfoldability Laver function, so that there is some
least $\theta$ with some $a=\alpha<\theta$ or $a=A\in H(\kappa^\plus)$ and a structure $N$ of size $\kappa$
with $\ell\in N$ such that there is no $\theta$-unfoldability embedding $j:N\to M$ with
$j(\ell)(\kappa)=a$. Fix an embedding $\bar j:V\to \bar M$ with critical point $\kappa$ and $\bar
j(\kappa)>\theta$. Since $P(\kappa)\of\bar M$, it follows that $\bar M$ agrees that there is such a
structure $N$ having no $\theta$-unfoldability embedding for which $\ell$ anticipates some $\alpha<\theta$
or $A\in H(\kappa^\plus)$. Since $\bar M$ may have fewer unfoldability embeddings than $V$, we know
$j(\ell)(\kappa)=a$ for some such choice of $N$ and $a<\theta'$ or $a\in H(\kappa^\plus)$, where
$\theta'\leq\theta$ and there is no $\theta'$-unfoldability embedding on $N$ for which $\ell$ anticipates
$a$. Consider the map $j=\bar j\restrict N$, so that $j:N\to \bar j(N)$. Note that $j(\ell)(\kappa)=\bar
j(\ell)(\kappa)=a$ and $j(\kappa)=\bar j(\kappa)>\theta\geq\theta'$. Thus, $j$ is a $\theta$-unfoldability
embedding on $N$ for which $\ell$ anticipates $a$. Furthermore, I claim that $j\in \bar M$. To see this,
code the structure $N$ with a relation $E$ on $\kappa$, so that $\<N,{\in}>\iso\<\kappa,E>$. Thus, $\<\bar
j(N),{\in}>\iso\<\bar j(\kappa),\bar j(E)>$. Furthermore, if $x\in N$ is represented by $\eta<\kappa$ with
respect to $E$, then $j(x)=\bar j(x)$ is represented by $\bar j(\eta)=\eta$ with respect to $\bar j(E)$.
That is, $j(x)=y$ if and only if $x$ has the same code below $\kappa$ with respect to $E$ that $y$ has with
respect to $\bar j(E)$. Thus, in the model $\bar M$ we can reconstruct the embedding $j:N\to \bar j(N)$
using $E$ and $\bar j(E)$. Thus, in $\bar M$ there is a $\theta'$-unfoldability embedding on $N$ for which
$\ell$ anticipates $a$ after all, a contradiction.\QED

One might want to improve this theorem to show that if $\kappa$ is unfoldable, then $\LDunf_\kappa$ holds,
by arguing in a similar manner; perhaps one would hope to use a big unfoldability embedding $\bar j:\bar
N\to \bar M$ in the place of $\bar j:V\to \bar M$ above. This strategy, however, seems to break down
because the model $N$ used in $\bar M$ to define $j(\ell)(\kappa)$ may not be in $\bar N$, and therefore we
won't be able to use the trick with $E$ and $\bar j(E)$ to know that the restricted embedding $j\restrict
N$ is in $\bar M$. What this idea does establish is that a weakened form of the Laver diamond holds for any
weakly compact or unfoldable cardinal (see Theorem \ref{WeakLDwc} below and the subsequent remarks). In a
forthcoming article, Mirna Dzamonja and I will show that $\LDunf_\kappa$ can fail for an unfoldable
cardinal $\kappa$, by showing that $\Diamond_\kappa(\Reg)$ can fail for such cardinals (extending Kai
Hauser's \cite{Hauser92:IndescribablesWithoutDiamond} result from $\Pi^m_n$-indescribable cardinals up to
strongly unfoldable cardinals).

One can improve Theorem \ref{MeasurableLDunf} by strengthening the conclusion still further:

\Theorem. If $\kappa$ is a measurable cardinal, then $\LDram_\kappa$ holds.

\Proof: Suppose that $\kappa$ is a measurable cardinal. Fix a well ordering $\trianglelt$ of $V_\kappa$ in
order type $\kappa$. I will define $\ell\from\kappa\to V_\kappa$ recursively. Suppose that
$\ell\restrict\gamma$ has been defined.  Let $\<a,b>$ be the $\trianglelt$-least pair, if any, such that
$b\of\gamma$ and there is no nice, iterable weakly amenable structure $\<M,F>$ of size $\gamma$ with
$\ell\restrict\gamma,b\in M$ and $a\in H(\kappa^\plus)\intersect M$ such that the corresponding ultrapower
map has $j(\ell\restrict\gamma)(\gamma)=a$. Let $\ell(\gamma)=a$, if such a pair exists.

I claim that $\ell$ witnesses $\LDram_\kappa$. If not, there is some set $B\of\kappa$ and some $A\in
H(\kappa^\plus)$ having no nice, iterable weakly amenable structure $\<M,F>$ containing $B$, $A$ and
$\ell$, whose ultrapower has $j(\ell)(\kappa)=A$. Let $\mu$ be any normal measure on $\kappa$, and consider
the corresponding ultrapower embedding $j:V\to M_0$. Since $M_0^\kappa\of M_0$, it follows that $M_0$ and
$V$ have the same structures of size $\kappa$, the same normal measures on such structures, and $M_0$
agrees with $V$ about whether such structures are iterable. In particular, the structure $M_0$ agrees with
$V$ that the pair $\<A,B>$ has no such structure anticipating $A$ via $j(\ell)\restrict\kappa=\ell$. Assume
without loss of generality this pair $\<A,B>$ is least with respect to $j(\trianglelt)$ with this property.
It follows that $j(\ell)(\kappa)=A$. For some suitably large $\theta$, choose any $X\elesub V_\theta$ of
size $\kappa$ with $X^\ltkappa\of X$ and $\ell,A,B,\mu\in X$. If $M$ is the Mostowski collapse of $X$, then
it follows that $F=\mu\intersect M\in M$ and $\<M,F>$ is an iterable weakly amenable structure of size
$\kappa$ with $\ell,A,B\in M$ and $M^\ltkappa\of M$ (iterable because $F$ is countably closed).
Furthermore, since the ultrapower of $X$ by $\mu$ has $j(\ell)(\kappa)=A$, it follows by the Mostowski
collapse that the ultrapower $j_0$ of $M$ by $F$ has $j_0(\ell)(\kappa)=A$, as desired.\QED

Let me now move on to automatic instances of $\LDmeas_\kappa$. A cardinal $\kappa$ is defined to be {\df
$\mu$-measurable} (see \cite{Mitchell1979:HypermeasurableCardinals}) when there is an embedding $j:V\to M$
with critical point $\kappa$ such that the induced normal measure $\mu=\set{X\of\kappa\st \kappa\in j(X)}$
is in $M$. Let me refer to such an embedding as {\df self-reflective}. Please note that a self-reflective
embedding $j$ cannot be the ultrapower by a normal measure, since by \cite[Lemma 28.9b]{Jech:SetTheory} no
measure is in its own ultrapower. But any $\kappa+2$-strong embedding $j:V\to M$, for example, is
self-reflective, since if $V_{\kappa+2}\of M$, then $M$ contains the same normal measures on $\kappa$ as
$V$ does, including the induced normal measure for $j$. Consequently, if $\kappa$ is $\kappa+2$-strong,
then it is $\mu$-measurable.

\Theorem. If\/ $\kappa$ is $\mu$-measurable, then $\LDmeas_\kappa$ holds.\label{SelfReflectiveMeasurable}

\Proof: Define $\ell$ as usual by transfinite recursion. Let $\ell(\gamma)$ be some $A\in H(\gamma^\plus)$
such that no ultrapower embedding $j:V\to M$ with critical point $\gamma$ has $j(\ell)(\gamma)=A$, if any
such $A$ exists. I claim that $\ell$ is a $\LDmeas_\kappa$ Laver function. If not, there is a set $A\in
H(\kappa^\plus)$ such that no ultrapower embedding $j:V\to M$ has $j(\ell)(\kappa)=A$. Using the hypothesis
that $\kappa$ is $\mu$-measurable, fix an embedding $j:V\to M$ with critical point $\kappa$ such that the
induced normal measure $\mu_0=\set{X\of\kappa\st\kappa\in j(X)}$ is in $M$. Since the measures on $\kappa$
in $M$ are all also measures in $V$, it follows that $M$ agrees that there is a set $A\in H(\kappa^\plus)$
which is not anticipated by $\ell=j(\ell)\restrict\kappa$ with respect to any ultrapower on $\kappa$. Thus,
$j(\ell)(\kappa)=A$ for some such set $A$. Let $j_0$ be the induced factor embedding by $\mu_0$, yielding
the commutative diagram \factordiagramup{V}{j}{M}{j_0}{k}{M_0} where $k$ has critical point above $\kappa$.
It follows that $j_0(\ell)(\kappa)=A$, and so the set $A$ is anticipated by $\ell$ via the ultrapower by
$\mu_0$. Since $M$ has $\mu_0$ and computes $j_0(\ell)$ the same as $V$ does, we see that $A$ is not a
counterexample after all, a contradiction.\QED

The key idea driving all the constructions of Laver functions up to now is that at each stage
$\ell(\gamma)$ is a minimal counterexample for the function $\ell\restrict\gamma$, and we have a larger
embedding $j:V\to M$ such that enough of the induced embedding $j_0:V\to M_0$, for whatever kind of
$\LDstar_\kappa$ is under discussion, is in $M$ so that it can compute $j_0(\ell)(\kappa)$ accurately. It
follows that $j(\ell)(\kappa)$ cannot be a counterexample, because any value $j(\ell)(\kappa)=A$ will
betray itself in the factored embedding $j_0(\ell)(\kappa)=A$. Thus, there can be no counterexample, and
the function $\ell$ will be a Laver function.

Using this idea, we can reduce the hypotheses on the local versions of supercompactness and strongness
above. Specifically, define that an elementary embedding $j:V\to M$ is {\df self-reflective} for $\theta$
supercompactness if and only if $j\image\theta\in M$ and the induced $\theta$-supercompactness measure
$\mu=\set{X\of P_\kappa\theta\st j\image\theta\in j(X)}$ is in $M$. (This is true, for example, if $j$ is a
$2^{\theta^\ltkappa}$-supercompactness embedding.) Similarly, the embedding $j:V\to M$ is {\df
self-reflective} for $\theta$-strongness if $V_\theta\of M$ and the induced $\theta$-strongness extender
$E=\set{\<X,s>\st X\of V_\kappa\And s\in j(X)\intersect \beth_\theta^\ltomega}$ is in $M$.

\Theorem. If $\kappa$ is $\theta$-supercompact and has an embedding that is self-reflective for
$\theta$-supercompactness, then $\LDthetasc_\kappa$. If it has an embedding that is self-reflective for
$\theta$-strongness, then $\LDthetastr_\kappa$.\label{SelfReflectiveSC}

It will follow from Theorem \ref{ForcingLDscLocal} that the converses of Theorems
\ref{SelfReflectiveMeasurable} and \ref{SelfReflectiveSC} needn't hold, since the consistency strength of
the Laver diamond is no greater than the existence of the corresponding large cardinal, whereas
self-reflective embeddings have a strictly greater consistency strength, since they imply the existence of
such large cardinals below the given cardinal. For example, the least measurable cardinal can never have a
self-reflective embedding, but it can have $\LDmeas_\kappa$.

It turns out that weakly compact cardinals can always exhibit a kind of self-reflectivity in their
embeddings (Lemma \ref{WCSelfReflection} below), and this can be used to prove a weakened form of the Laver
diamond for every weakly compact cardinal. The resulting function $\ell$ still anticipates every set $A\in
H(\kappa^\plus)$, but one has somewhat less freedom to choose the domain of the embedding. Specifically,
while $\LDwc_\kappa$ makes an assertion with the quantifier structure $\exists\ell\,\forall A\,\forall
M\,\exists j\, (j(\ell)(\kappa)=A)$, the following theorem establishes that every weakly compact cardinal
exhibits at least $\exists\ell\,\forall A\,\exists M\,\exists j\,(j(\ell)(\kappa)=A)$.

\Theorem. If $\kappa$ is a weakly compact cardinal, then for any set $B\of\kappa$ there is a function
$\ell\from\kappa\to V_\kappa$ such that for any set $A\in H(\kappa^\plus)$ there is a nice structure $M$
with $B,\ell\in M$ and an embedding $j:M\to N$ such that $j(\ell)(\kappa)=A$.

\Proof: The first step is to prove that weakly compact cardinals exhibit a kind of self-reflectivity.

\SubLemma.(Hauser) A cardinal $\kappa$ is weakly compact if and only if for any transitive structure $M$ of
size $\kappa$ with $M^\ltkappa\of M$, there is an embedding $j:M\to N$ into a transitive set $N$, with
critical point $\kappa$, such that $M$ and $j$ are elements of $N$.\label{WCSelfReflection}

\Proof: By insisting that $M$ and $j$ are in $N$, we have strengthened the usual embedding characterization
of weakly compact cardinals, so I need only prove the forward direction. Suppose that $\kappa$ is weakly
compact, and that $M$ is a transitive structure of size $\kappa$ with $M^\ltkappa\of M$. Let $\bar M$ be a
transitive structure of size $\kappa$ with $M\in\bar M$ and $\bar M\satisfies|M|=\kappa$. There is
therefore some relation $E$ on $\kappa$ such that $\<M,{\in}>\iso\<\kappa,E>$. This isomorphism is unique,
since it is the Mostowski collapse of $\kappa$ under $E$. Fix an embedding $j:\bar M\to \bar N$ with
critical point $\kappa$, and let $j_0=j\restrict M$ and $N=j(M)$, so that $j_0:M\to N$. Since
$E=j(E)\intersect\kappa\cross\kappa$, it follows that $E\in\bar N$, and consequently also $M\in \bar N$. Of
course, $j(E)\in\bar N$ as well. Furthermore, if $a\in M$ is coded by $\alpha$ with respect to $E$, then
$j(a)$ is coded by $\alpha$ ($=j(\alpha)$) with respect to $j(E)$. Therefore, we can reconstruct $j_0$ from
$E$ and $j(E)$ in $\bar N$: simply map the set coded by any $\alpha<\kappa$ with respect to $E$ to the set
coded by $\alpha$ with respect to $j(E)$. Thus, $j_0\in\bar N$. Finally, since $M^\ltkappa\of M$, and this
is true in $M$, we know that $\bar N\satisfies j(M)^{{<}j(\kappa)}\of j(M)$. In particular, $N=j(M)$ has
the same structures of size $\kappa$ as $\bar N$, and so $M\in N$. This implies in addition that $M\of N$
and since $j_0:M\to N$, we actually know that $j_0\of N$, considered as a set of ordered pairs. Since this
embedding has size $\kappa$ (the size of its domain), it follows that $j_0\in N$. So we have $j_0:M\to N$
with $j_0\in N$, as desired.

Please note that $\union\ran(j_0)$ is constructible from $j_0$ in $N$. Therefore, the embedding $j_0$
cannot be cofinal in $N$. Nevertheless, the embedding $j_0:M\to\union\ran(j_0)$ into the truncated target
is in $N$, and fully elementary.\QED

To prove the theorem, I define $\ell$ recursively as before. Fix $B\of \kappa$ and suppose that
$\ell\restrict\gamma$ is defined. Let $\ell(\gamma)$ be any set $A\in H(\gamma^\plus)$ such that there is
no transitive structure $M$ of size $\gamma$ with $\ell\restrict\gamma$ and $B\intersect\gamma$ in $M$ and
an embedding $j:M\to N$ with critical point $\gamma$ such that $j(\ell\restrict\gamma)(\gamma)=A$, if any
such set $A$ exists. I claim that this function $\ell\from\kappa\to V_\kappa$ satisfies the conclusion of
the theorem. If not, there is some $A\in H(\kappa^\plus)$ not having any structure and embedding that
anticipates it. Fix any structure $M$ containing $B$, $A$ and $\ell$, and by the Lemma fix an embedding
$j:M\to N$ with $M$ and $j$ in $N$. Since $N$ is transitive and $\ell,A\in N$, it also can have no
structure and embedding $j_0:M_0\to N_0$ with $B\in M_0$ and $j_0(\ell)(\kappa)=A$, since there is no such
structure and embedding in $V$. Therefore, $\kappa$ is in the domain of $j(\ell)$. Let
$A'=j(\ell)(\kappa)$. By the definition of $\ell$, this means that in $N$ there is no structure $M_0$
containing $j(B)\intersect\kappa$ ($=B$) and $j(\ell)\restrict\kappa$ ($=\ell$) and embedding $j_0:M_0\to
N_0$ with $j_0(\ell)(\kappa)=A'$. But $j:M\to \union\ran(j)$ is exactly such an embedding in $N$, a
contradiction.\QED

The point of the set $B$ in the theorem is to ensure, by insisting $B\in M$, that the structure $M$ has a
certain minimal amount of desired information. But actually, since the set $B$ can be coded directly into
the function $\ell$, and then from $\ell\in M$ one may deduce $B\in M$ as well, the version of the theorem
without $B$ easily implies the version with $B$ stated above.

The degree of control over $M$ provided by the theorem is not sufficient to carry out the usual proof of
$\Diamond_\kappa(\Reg)$ as a consequence, because in that argument one fixes a set $A\of \kappa$ and a club
$C$ on which $\ell$ does not anticipate it, and then gets a contradiction if $C\in M$ because
$j(\ell)(\kappa)=A=j(A)\intersect\kappa$ and $\kappa\in j(C)$.  This argument would not work if you had to
fix the club $C$ before knowing $\ell$ or $A$. Indeed, it is known that $\Diamond_\kappa(\Reg)$ can fail at
weakly compact, indescribable, unfoldable and strongly unfoldable cardinals (see Theorem \ref{LDwcCanFail}
below and the subsequent remarks), whereas the above theorem establishes the weak Laver diamond at every
weakly compact cardinal.

A class version of this argument establishes that there is a class function $\ell\from\ORD\to V$ such that
for any weakly compact cardinal $\kappa$ and any $A\in H(\kappa^\plus)$ there is a structure $M$ with
$\ell\restrict\kappa\in M$ and an embedding $j:M\to N$ with $j(\ell\restrict\kappa)(\kappa)=A$.
Furthermore, this construction works with unfoldable, indescribable and strongly unfoldable cardinals just
as easily (including the analogue of Lemma \ref{WCSelfReflection} for those cardinal notions), and one can
obtain the analogous weakening of the Laver diamond for these cardinals.

Next, I point out an obvious connections between stronger and weaker diamond principles.

\Observation. If $*_1$ and $*_2$ are two large cardinal notions and $*_1$ is at least as strong as $*_2$ in
the sense that every $*_1$ embedding is also a $*_2$ embedding, and the class of sets to be anticipated by
a $*_1$ Laver function includes all the sets that must be anticipated by a $*_2$ Laver function, then every
$*_1$ Laver function is also a $*_2$ Laver function.\label{TrivialObservation}

Because of this, we know that every $\LDsc_\kappa$ Laver function, for example, is a $\LDstr_\kappa$ Laver
function, every $\LD_\kappa^{\hbox{\tiny\!\!$\beth_\theta$-sc}}$ Laver function is a $\LDthetastr_\kappa$
Laver function, and so on. One can also show that every $\LDmeas_\kappa$ Laver function is a
$\LDunf_\kappa$ Laver function, by iterating an appropriate normal measure on $\kappa$, so as to send
$j(\kappa)$ above any desired $\theta$.

Let me now consider some seemingly weaker forms of the Laver diamond, analogous to the classical
$\Diamond^-$, where one has a small collection of sets at each stage rather than a single set. The fact is
that these apparently weaker forms are fully equivalent to the Laver diamond principle.

\Theorem. The following are equivalent.
\begin{enumerate}
\item $\LDwc_\kappa$.
\item There is a function  $L\from\kappa\to V_\kappa$ such that $|L(\gamma)|\leq\gamma$ for all
$\gamma\in\dom(L)$, and for any set $A\in H(\kappa^\plus)$ and any nice structure $M$ of size $\kappa$
containing $L$ and $A$, there is an embedding $j:M\to N$, having critical point $\kappa$, with $A\in
j(L)(\kappa)$.
\item There is a family $\cal L$ of at most $\kappa$ many functions from $\kappa$ to $V_\kappa$ such that
for any set $A\in H(\kappa^\plus)$ and any nice structure $M$ of size $\kappa$ containing $A$ and $\cal L$,
there is an embedding $j:M\to N$, having critical point $\kappa$, such that $j(\ell)(\kappa)=A$ for some
$\ell\in{\cal L}$.
\end{enumerate}\label{WeakLDwc}

\Proof: This argument adapts the classical proof that $\Diamond^-$ is equivalent to $\Diamond$ (see, for
example, \cite[Theorem 7.14]{Kunen:Independence}). Clearly statement 1 implies statement 2. Let me argue
that statement 2 implies statement 3. Fix $L$ as in statement 2, and for each $\gamma\in\dom(L)$, fix an
enumeration of $L(\gamma)$ in order type at most $\gamma$. Let $\ell_\beta(\gamma)$ be the $\beta^\th$
element of $L(\gamma)$, if it exists, and ${\cal L}=\set{\ell_\beta\st\beta<\kappa}$. I claim that this
family of functions witnesses statement 3. To see this, suppose that $A\in H(\kappa^\plus)$ and $M$ is a
nice structure of size $\kappa$ containing $A$ and $\cal L$. It follows that $L\in M$ as well, and the
orderings of each $L(\gamma)$ are constructible in $M$ from $\cal L$. Consequently, by statement 2, there
is an embedding $j:M\to N$ with $A\in j(L)(\kappa)$. Thus, $A$ is the $\beta^\th$ element of $j(L)(\kappa)$
for some $\beta<\kappa$, and therefore $A=j(\ell_\beta)(\kappa)$, as desired.

Let me now prove that statement 3 implies statement 1. Suppose that ${\cal
L}=\set{\ell_\beta\st\beta<\kappa}$ is a family of functions witnessing statement 3. Define the function
$\ell^*_\beta(\gamma)=x$ if $\ell_\beta(\gamma)$ is an $\gamma$-sequence whose $\beta^\th$ element is $x$.
I claim that one of the functions $\ell^*_\beta$ is a Laver function witnessing $\LDwc_\kappa$. If not,
then for each $\beta<\kappa$ there is a set $A_\beta$ and a structure $M_\beta$ having no embedding
$j:M_\beta\to N$ for which $j(\ell^*_\beta)(\kappa)=A_\beta$. Let $A=\<A_\beta\st\beta<\kappa>$ and choose
some nice structure $M$ containing $A$, $\cal L$ and all the structures $M_\beta$. By statement 3, there is
an embedding $j:M\to N$ with $j(\ell_\beta)(\kappa)=A$ for some $\beta<\kappa$. Since $A$ is a
$\kappa$-sequence whose $\beta^\th$ element is $A_\beta$, this means that
$j(\ell^*_\beta)(\kappa)=A_\beta$. By restricting the embedding $j$ to $M_\beta$, we conclude that
$A_\beta$ and $M_\beta$ were anticipated by $\ell^*_\beta$ after all, a contradiction.\QED

\Theorem. The weak forms of $\LDwc_\kappa$, $\LD^{\hbox{\!\!\tiny$\Pi^m_n$-ind}}_\kappa$, $\LDunf_\kappa$,
$\LDsunf_\kappa$, $\LDram_\kappa$, $\LDmeas_\kappa$, $\LDstr_\kappa$, $\LDsuperstrong_\kappa$ and
$\LDsc_\kappa$, analogous to those stated in Theorem \ref{WeakLDwc}, namely, where $L(\gamma)$ is a
collection of $\gamma$ many sets to be anticipated or alternatively, where there are $\kappa$ many
proto-Laver functions $\ell$ such that any set $A$ is anticipated by at least one of them, are equivalent
to the full versions of the principles, respectively.\label{WeakLD}

\Proof: The argument of Theorem \ref{WeakLDwc} generalizes easily to the other large cardinals. What is
needed is that the class of sets $A$ to be guessed must be closed under $\kappa$-sequences (or at least
closed under codes for such sequences), so that the counterexamples $A_\beta$ for $\ell^*_\beta$ can be
combined into the sequence $A=\<A_\beta\st\beta<\kappa>$. Then, for each of the large cardinal embedding
types, if $j(\ell_\beta)(\kappa)=A$, one directly concludes $j(\ell^*_\beta)(\kappa)=A_\beta$,
contradicting the choice of $A_\beta$.\QED

If one allows $L(\gamma)$ to guess more than $\gamma$ many sets, the resulting weak Laver principle can
become trivial. For example, the function $L(\gamma)=H(\gamma^\plus)$ clearly allows one to capture any set
$A\in H(\kappa^\plus)$. When the \GCH\ fails, one might consider guessing with sets of size $\gamma^\plus$
or more, and the exact strength of such a principle is not clear.

Let me conclude this section by mentioning a stronger form of the Laver diamond principle, essentially
appearing in \cite{Hamkins2000:LotteryPreparation}. Let $\LDLD_\kappa$ be the principle that asserts that
there is a function $\ell\from\kappa\to V_\kappa$ such that for any (appropriate) embedding $j:V\to M$ and
any set $A\in M_{j(\kappa)}$, there is another embedding $h:V\to M$ (the same $M$) with $h(\ell)(\kappa)=A$
and $h\restrict\ORD=j\restrict\ORD$. The class of appropriate embeddings $j$ is determined by whatever
large cardinal context one has at hand, leading to principles $\LDLD_\kappa^{\hbox{\!\!\tiny meas}}$ and so
on. This principle makes sense in a variety of large cardinal contexts, as well as for set-sized embeddings
$j:M\to N$. The results of \cite{Hamkins2000:LotteryPreparation}, mentioned in Theorems \ref{Fast Function
Flexibility Theorem} and \ref{Generalized Laver Function Theorem} below, implicitly concern the possibility
of forcing such a principle.

\Section Laver functions in $L$ and $L$-like models\label{LDinL}

For the most part, large cardinals in the constructible universe $L$ have Laver functions there. Let me
begin with the weakly compact cardinals. I am indebted to James Cummings and Ernest Schimmerling for the
following observation, arising in an email exchange about these notions. An equivalent result appeared
meanwhile, independently, in \cite{ShelahVaananan726:ExtensionsOfInfinitaryLogic}.

\Theorem. If\/ $V=L$, then $\LDwc_\kappa$ holds for every weakly compact cardinal $\kappa$. Indeed, $\LDwc$
holds: there is a universal Laver function $\ell\from\ORD\to L$ such that $\ell\restrict\kappa$ is a weak
compactness Laver function for any weakly compact cardinal $\kappa$.\label{LDwcinL}

\Proof: Again, we define $\ell$ by transfinite recursion. To define $\ell(\gamma)$, let
$\<M_\gamma,A_\gamma>$ be the $L$-least pair, if any, such that $M_\gamma$ is a transitive set of size
$\gamma$ containing $\gamma$, $A_\gamma$ and $\ell\restrict\gamma$ and such that $A_\gamma\in
H(\gamma^\plus)$ is not anticipated by $\ell\restrict\gamma$ for any weak compactness embedding
$j:M_\gamma\to N$ with critical point $\gamma$. That is, $j(\ell\restrict\gamma)(\gamma)\not=A_\gamma$ for
all such embeddings $j$. Define $\ell(\gamma)=A_\gamma$. I claim now that $\ell$ is a weak compactness
Laver function. If not, there is some $L$-least pair $\<M,A>$ with $\kappa,A,\ell\in M$ such that no weak
compactness embedding $j:M\to N$ has $j(\ell)(\kappa)=A$. Fix a transitive set $\bar M$ of size $\kappa$
with $M\in \bar M$, such that $\bar M$ agrees that $\<M,A>$ is the $L$-least such pair. (To do this, simply
take an $\bar M$ that has embeddings $h:M'\to N'$ with $h(\ell)(\kappa)=A'$ for all $L$-smaller pairs
$\<M',A'>$, of which there are only $\kappa$ many.) Now let $j:\bar M\to \bar N$ be any weak compactness
embedding on $\bar M$. Since $P(\kappa)^{\bar M}\of \bar N$, it follows that $\bar N$ also has all those
embeddings $h:M'\to N'$, and so $\bar N$ agrees that $\<M,A>$ is the least such for which there is no
embedding. Thus, $j(\ell)(\kappa)=A$. Since the restricted embedding $j\restrict M:M\to j(M)$ still has
$j(\ell)(\kappa)=A$, this contradicts the choice of the pair $\<M,A>$. So $\ell$ witnesses $\LDwc_\kappa$.

For the second sentence of the theorem, simply note that the definition of $\ell(\gamma)$ did not depend on
$\kappa$. So the Laver functions I constructed for the various weakly compact cardinals all cohere into a
universal Laver function.\QED

\Theorem. If\/ $V=L$, then $\LDunf_\kappa$ holds for every unfoldable cardinal $\kappa$, and there is a
universal $\LDunf$ Laver function. Indeed, these are ordinal-anticipating Laver functions, and they witness
$\LDsunf$ as well.\label{LDsunfinL}

\Proof: We follow the construction of Theorem \ref{LDwcinL}. If $\ell$ is defined up to $\gamma$, let
$\<M,\theta,a>$ be the least triple, if any, such that $M$ is a transitive set of size $\gamma$ with
$\ell\restrict\gamma\in M$, the set $a$ is in $L_\theta$, and there is no $\theta$-unfoldability embedding
$h:M\to N$ such that $h(\ell)(\gamma)=a$. Define $\ell(\gamma)=a$. A simple reflection argument shows that
$\theta$ is below the next strongly unfoldable cardinal (which, in $L$, is the same as the next unfoldable
cardinal). I claim that for any unfoldable cardinal $\kappa$ the resulting function $\ell\from\kappa\to
V_\kappa$ is a $\LDunf_\kappa$ Laver function, anticipating every element of $L_\theta$, for any ordinal
$\theta$, with a $\theta$-unfoldability embedding. If not, there is some least triple $\<M,\theta,a>$
forming a counterexample. Fix any transitive $\bar M$ of size $\kappa$ with $M\in \bar M$, and select
$\lambda>\theta$ such that all the witnessing embeddings for earlier triples $\<M',\theta',a'>$ appear
before $\lambda$ in the canonical ordering of the universe. Since $\kappa$ is unfoldable, there is an
embedding $j:\bar M\to \bar N$ with critical point $\kappa$ and $j(\kappa)>\lambda$. By the choice of
$\lambda$, we know that $\<M,\theta,a>\in\bar N$, and furthermore, the structure $\bar N$ knows that no
earlier triple is a counterexample. Thus, since the triple $\<M,\theta,a>$ remains a counterexample in
$\bar N$, it follows that $j(\ell)(\kappa)=a$. Let $h=j\restrict M$, so that $h:M\to j(M)$. Since
$h(\kappa)=j(\kappa)>\lambda\geq\theta$ and $h(\ell)(\kappa)=j(\ell)(\kappa)=a$, we know that $h$ is a
$\theta$-unfoldability embedding on $M$ for which $\ell$ anticipates $a$, contradicting our choice of
$\<M,\theta,a>$.

The final claim of the theorem is proved with the observation that the notions of strong unfoldability and
unfoldability coincide in $L$, because including $V_\theta$ into the target of an embedding amounts to
including $L_{\beth_\theta}$, which amounts to the target simply including the ordinals up to
$\beth_\theta$. Since the Laver functions we just constructed anticipate every element of
$L_{\beth_\theta}$ with $\beth_\theta$-unfoldability embeddings, they will witness $\LDsunf$.\QED

These results generalize out of the constructible universe. The key idea at play is that there is a
definable well-ordering of $L$ that is absolute to the various models in which it is consulted. So let me
define that a well order $\trianglelt$ of the universe is {\df absolutely definable} by the formula
$\varphi$ if $\varphi$ defines $\trianglelt$ and any transitive (set) model $M$ in which $\varphi$ defines
a well order agrees with $V$ on whether $a\trianglelt b$, for any two objects $a,b\in M$. The canonical
order of $L$, for example, is absolutely definable by $a\trianglelt b\Iff (V=L$ and $a$ is constructed
before $b$), since any two transitive models of $V=L$ agree on the order in which elements of $L$ are
constructed. The restriction to set models makes the notion transparently first order, but it can be
omitted, since whenever $N$ is a proper class in which $\varphi$ defines a well order, then we obtain
absoluteness to $N$ by applying the absoluteness to sufficiently reflective set-sized substructures
$M\elesub_n N$. The theorems below can therefore be viewed as theorem schemes, ranging over the possible
definitions of $\trianglelt$.

\Theorem. Suppose that $\trianglelt$ is an absolutely definable well-ordering of the universe.  Then there
is a universal $\LDwc$ Laver function. If the definition of $\trianglelt$ uses a parameter, then the
universal $\LDwc$ function works for weakly compact cardinals above the parameter.\label{WOLDwc}

\Proof: Let us try to proceed just as in Theorem \ref{LDwcinL}. Suppose that $\trianglelt$ is an absolutely
definable well-order of $V$ (defined with parameter $z$) and that $\ell$ has been defined up to $\gamma$.
Let $\<M_\gamma,A_\gamma>$ be the $\trianglelt$-least pair, if any, such that $M_\gamma$ is a transitive
set of size $\gamma$ with $\gamma,A_\gamma,\ell\restrict\gamma\in M_\gamma$ and $A_\gamma\in
H(\gamma^\plus)$ is not anticipated by $\ell\restrict\gamma$ for any weak compactness embedding
$j:M_\gamma\to N$ with critical point $\gamma$. Define $\ell(\gamma)=A_\gamma$. Once again, I claim that
$\ell$ is a universal $\LDwc$ Laver function. If not, there is some weakly compact cardinal $\kappa$ (with
the parameter $z$ in $V_\kappa$) and a $\trianglelt$-least pair $\<M,A>$ with $\kappa,A,\ell\in M$ such
that no weak compactness embedding $j:M\to N$ has $j(\ell\restrict\kappa)(\kappa)=A$. Fix a transitive set
$\bar M$ of size $\kappa$ with $M\in \bar M$, such that $\bar M$ agrees that $\<M,A>$ is the
$\trianglelt$-least counterexample (this can be done by collapsing a suitable elementary substructure of
some large $V_\theta$). Thus, $\bar M$ has the embeddings for all the $\trianglelt$-earlier pairs,
witnessing that they are anticipated by $\ell\restrict\kappa$. Since these embeddings are coded by subsets
of size $\kappa$, it follows that they are also in $\bar N$. Since the pair $\<M,A>$ remains a
counterexample in $\bar N$, and by the absoluteness of $\trianglelt$ the predecessors of this object in
$\bar N$ are among the actual predecessors in $V$, for all of which $\bar N$ has the witnessing embeddings,
it follows that $\bar N$ also agrees that $\<M,A>$ is the least pair for which there is no embedding. Thus,
$j(\ell\restrict\kappa)(\kappa)=A$. Since the restricted embedding $j\restrict M:M\to j(M)$ still has
$j(\ell\restrict\kappa)(\kappa)=A$, this contradicts the choice of the pair $\<M,A>$. So
$\ell\restrict\kappa$ witnesses $\LDwc_\kappa$, and the proof is complete.\QED

The argument also generalizes in the case of strongly unfoldable cardinals.

\Theorem. Suppose that $\trianglelt$ is an absolutely definable well-ordering of the universe. Then there
is a universal $\LDsunf$ Laver function. If the definition uses a parameter, then the universal $\LDsunf$
Laver function works for strongly unfoldable cardinals above the parameter.\label{AbsWOLDsunf}

\Proof: We follow Theorem \ref{LDsunfinL}. Suppose that $\trianglelt$ is an absolutely definable
well-ordering of the universe (with parameter $z$) and that $\ell$ has been defined up to $\gamma$. Let
$\<M_\gamma,\theta_\gamma,A_\gamma>$ be the $\trianglelt$-least triple, if any, such that $M_\gamma$ is a
transitive set of size $\gamma$ with $\ell\restrict\gamma\in M_\gamma$ and $A_\gamma\in
H(\theta_\gamma^\plus)$ and there is no $\theta_\gamma$-strong unfoldability embedding $h:M_\gamma\to N$
such that $h(\ell)(\gamma)=A_\gamma$. Let $\ell(\gamma)=A_\gamma$. If $\kappa$ is a strongly unfoldable
cardinal and the parameter $z$ is in $V_\kappa$, then I claim that $\ell^*=\ell\restrict\kappa$ is a
$\LDsunf_\kappa$ Laver function. If not, there is some least triple $\<M,\theta,A>$ forming a
counterexample. Fix any transitive $\bar M$ of size $\kappa$ with $M\in \bar M$, and select
$\lambda>\theta$ such that all the witnessing embeddings for earlier triples $\<M',\theta',A'>$ appear in
$H(\lambda^\plus)$. Since $\kappa$ is strongly unfoldable, there is an embedding $j:\bar M\to \bar N$, with
critical point $\kappa$, such that $j(\kappa)>\lambda$ and $H(\lambda^\plus)\of\bar N$. By the choice of
$\lambda$, we know that $\<M,\theta,A>\in\bar N$, and furthermore, $\bar N$ agrees with $V$ on
$\trianglelt\restrict H(\theta^\plus)$. Since all the embeddings for the earlier triples appear in $\bar
N$, this model therefore knows that no earlier triple is a counterexample, and since $\<M,\theta,A>$
remains a counterexample in $\bar N$, it follows that $j(\ell)(\kappa)=A$. Let $h=j\restrict M$, so that
$h:M\to j(M)$. Since $h(\kappa)=j(\kappa)>\lambda\geq\theta$ and $h(\ell)(\kappa)=j(\ell)(\kappa)=A$, we
know that $h$ is a $\theta$-unfoldability embedding on $M$ for which $\ell$ anticipates $A$. Furthermore,
since $M$ and $\bar M$ agree on $H(\kappa)$, it follows that $\bar N$ and $j(M)$ agree on $H(j(\kappa))$,
so in particular $H(\theta^\plus)\of j(M)$; thus, $h$ is a $\theta$-strong unfoldability embedding. So
there is a $\theta$-unfoldability embedding on $M$ for which $\ell$ anticipates $A$ after all, a
contradiction.\QED

This result can be improved by relaxing the hypotheses on the order somewhat. It is not necessary that
$\trianglelt$ is absolute to {\it all} transitive models in which the definition defines a well-order, but
rather, the order only needs to agree on $H(\theta^\plus)$ for models that have the correct
$H(\theta^\plus)$. For example, if the order $\trianglelt$ were defined locally, in the sense that
$a\trianglelt b\Iff H(\theta^\plus)\satisfies\varphi(a,b)$, where $\theta=\max\{|\TC(a)|,|\TC(b)|\}$, then
this would be the case. More generally, let me define that a definable well-order $\trianglelt$ is {\df
locally absolute} if for any cardinal $\theta$, any transitive (set) model with the correct
$H(\theta^\plus)$ in which the definition defines a well-order agrees on $\trianglelt\restrict
H(\theta^\plus)$.

\Theorem. Suppose that $\trianglelt$ is a definable well-ordering of the universe that is locally absolute.
Then there is a universal $\LDsunf$ Laver function. If the definition uses a parameter, then the universal
$\LDsunf$ Laver function works for strongly unfoldable cardinals above the
parameter.\label{LocallyAbsWOLDsunf}

\Proof: A careful scrutiny will reveal that this is all that is used in the proof of Theorem
\ref{AbsWOLDsunf}.\QED

One shouldn't take away from the results in this section the expectation that Laver functions always exist
in canonical inner models for large cardinals. For example, Theorem \ref{LDmeasFailsInL[mu]} will show that
$\LDmeas_\kappa$ fails in $L[\mu]$, the canonical inner model of a measurable cardinal. The Laver diamond
$\LDthetastr_\kappa$ similarly fails in the canonical model for a $\theta$-strong cardinal $\kappa$ (e.g.
$(\kappa+3)$-strong but not $(\kappa+4)$-strong), but in the canonical model for a fully strong cardinal,
$\LDstr$ holds of course by Corollary \ref{UniformLDstr}, since strong cardinals always have Laver
functions.

Let me now consider the ordinal-anticipating Laver functions. Recall that these are the functions which are
required to anticipate only the ordinals (up to some bound $\delta$), with no assumption that they
anticipate other kinds of sets; thus, they may not be Laver diamond functions at all. For example, in the
case of measurability a function $f\from\kappa\to V_\kappa$ (or more to the point, $f\from\kappa\to\kappa$;
but I want to allow the possibility that $f$ is a full Laver diamond function) is an {\df ordinal
anticipating} Laver function up to $\delta$ if for every ordinal $\alpha<\delta$ there is an embedding
$j:V\to M$ with critical point $\kappa$ such that $j(f)(\kappa)=\alpha$. Natural choices for $\delta$ would
include $2^\kappa$, $\kappa^\plus$, $\kappa$, or even $\ORD$ itself, and it is not entirely clear yet how
all these hypotheses are related. There are natural analogues of this concept for the other large cardinal
notions.

\Question. If there is an ordinal-anticipating Laver function, does the corresponding Laver diamond
principle hold? For example, if $\kappa$ is measurable and there is a function anticipating all ordinals up
to $2^\kappa$, then does $\LDmeas_\kappa$ hold? If $\kappa$ is unfoldable and there is an ordinal
anticipating Laver function for unfoldability embeddings, must $\LDunf_\kappa$ hold?
\label{OrdinalAnticipatingQuestion}

One has such questions for any of the Laver diamond principles. If there is a sufficiently absolute well
ordering of the universe, then the theorems below provide an affirmative answer. Later, in section
\ref{ForcingLD}, I will prove that when one forces to add an ordinal-anticipating function, a full
set-anticipating Laver function can be constructed from it. Proving an outright implication would clarify
matters.

\Theorem. Suppose that $\trianglelt$ is a definable well ordering of the universe that is locally absolute.
If $\kappa$ is a measurable cardinal and has a ordinal anticipating Laver function, anticipating all
ordinals up to the order type of $H(\kappa^\plus)$ with respect to $\trianglelt$, then $\LDmeas_\kappa$
holds.\label{OrdinalOnlyAnticipatingLDmeas}

\Proof: I will actually need only that $\trianglelt$ is a definable well ordering of the universe (using
parameters in $V_\kappa$), whose ordering of $H(\kappa^\plus)$ is absolute to models with the correct
$H(\kappa^\plus)$. Suppose that $f\from\kappa\to\kappa$ anticipates all ordinals up to the order type
$\delta$ that $H(\kappa^\plus)$ forms as a subset of the field of $\trianglelt$. To define the full Laver
function, simply let $\ell(\gamma)$ be the $f(\gamma)^\th$ element of $H(\gamma^\plus)$ in the
$\trianglelt$ order. To see that this works, consider any set $A\in H(\kappa^\plus)$. This set is the
$\alpha^\th$ element of $H(\kappa^\plus)$ with respect to the order $\trianglelt$ for some $\alpha<\delta$.
Let $j:V\to M$ be any embedding with critical point $\kappa$ such that $j(f)(\kappa)=\alpha$. Since $M$ and
$V$ have the same $H(\kappa^\plus)$ and $j(\trianglelt)$ is defined in $M$ by the same formula as
$\trianglelt$ in $V$, it follows by the local absoluteness hypothesis on $\trianglelt$ that the
$j(\trianglelt)$ order on $H(\kappa^\plus)$ is the same as the $\trianglelt$ order in $V$. Thus,
$j(\ell)(\kappa)=A$, as desired.\QED

I do not know at the moment whether the definable well ordering hypothesis can be omitted. It would be very
interesting to have a model with a function that could anticipate any ordinal, but no set anticipating
Laver function. Theorem \ref{OrdinalOnlyAnticipatingLDmeas} clearly generalizes to the Laver diamond
principles in many other large cardinal contexts.

There is also a kind of converse question. If $\kappa$ is a measurable cardinal, for example, then there
are embeddings sending $\kappa$ arbitrarily high in the ordinals, and so it is conceivable that there could
be a function anticipating all ordinals. This idea should be tempered by the realization that embeddings
arising as iterations of a normal measure never help you anticipate more ordinals, since all iterations
agree on $j(f)(\kappa)$ with the original one-step ultrapower. But the question remains:

\Question. If\/ $\LDmeas_\kappa$ holds, is there an ordinal anticipating Laver function
$f\from\kappa\to\kappa$ for measurability that anticipates all ordinals?

One might alternatively like to assume only that there is a function anticipating all ordinals up to some
$\delta$, and ask whether this bound can be freely removed. Of course, by reflection there is some bound
$\delta$ making an affirmative answer, but what we really would like to know is whether the question is
true for some natural choices of $\delta$, such as $2^\kappa$.

\Observation. A particular function might anticipate all ordinals $\alpha<2^\kappa$, but not all ordinals.

\Proof: Suppose that $f$ anticipates all $\alpha<2^\kappa$, and let $f^*$ be the restriction of $f$ to the
set of $\gamma$ for which $f(\gamma)<2^\gamma$. For any $\alpha<2^\kappa$ there is $j:V\to M$ with
$j(f)(\kappa)=\alpha$, and since $M$ agrees that $\alpha<2^\kappa$, it follows that $\kappa$ is in the
domain of $j(f^*)$, and so $j(f^*)(\kappa)=j(f)(\kappa)=\alpha$. Thus, $f^*$ anticipates every
$\alpha<2^\kappa$. But $f^*$ cannot anticipate all ordinals, because if there is an embedding $j:V\to M$
with $j(f^*)(\kappa)=\beta$, then by the definition of $f^*$, we know that $\beta<(2^\kappa)^M$, and this
bound must be less than $(2^\kappa)^\plus$.\QED

This cutting-off phenomenon occurs with larger (definable) bounds just as well, but it doesn't seem to get
to the heart of the question of whether there could be other functions that do anticipate all ordinals.
What one wants to know is whether there could be a function that anticipates all ordinals up to $2^\kappa$,
when no function anticipates all ordinals.

For strongly unfoldable cardinals, it seems natural to ask the following question.

\Question. If\/ $\kappa$ is strongly unfoldable and $\LDunf_\kappa$, with ordinal anticipating Laver
functions, then must $\LDsunf_\kappa$ hold?

If $V=L$, then the answer is trivially yes, because $\LDunf_\kappa$ and $\LDsunf_\kappa$ are equivalent in
$L$ (and both hold there). In the general case, one might hope to start with an ordinal-anticipating
$\LDunf_\kappa$ Laver function $\ell\from\kappa\to\kappa$ and construct from it a function
$\ell^*\from\kappa\to V_\kappa$ by defining $\ell^*(\gamma)$ to be the $\ell(\gamma)^\th$ element of
$V_\kappa$, using some well-order $\trianglelt$ of $V_\kappa$. The difficulty with this approach, however,
is that when using it one needs to control both $j(\ell)(\kappa)$ and $j(\trianglelt)$ in order to ensure
that $j(\ell^*)(\kappa)$ is as desired. Of course, if there is a locally absolute definable well-ordering
of the universe, then one can simply use that order, and the idea will work. But by Theorem
\ref{LocallyAbsWOLDsunf} we already know $\LDsunf$ in this case directly, without needing to assume an
ordinal-anticipating $\LDunf_\kappa$ function. Can one get by without the definable well-ordering of the
universe?

\Section Forcing the Laver diamond principles $\LaverDiamond_\kappa$\label{ForcingLD}

For a great variety of large cardinals---weakly compact, indescribable, unfoldable, strongly unfoldable,
Ramsey and measurable cardinals, for example---one can force the existence of a Laver function.  And unlike
the classical results and the results of the previous section, where we used a greater large cardinal
property to deduce the existence of a Laver function for a lesser large cardinal property, these forcing
arguments generally allow one to add a Laver function without sacrificing any degree of large cardinal
strength. Therefore, in most cases the principle $\LDstar_\kappa$ is equiconsistent with the large cardinal
property $*$ itself.

For each of the large cardinal notions, the forcing to add a Laver function is the same: Woodin's fast
function forcing. (Later, I will show how to force the Laver diamond principle via Silver forcing and even
via $\kappa$-c.c. forcing.) Please consult \cite{Hamkins2000:LotteryPreparation} for a careful introduction
to fast function forcing; here, I will quickly develop the basic facts. Suppose $\kappa$ is an inaccessible
cardinal. Conditions in the {\df fast function forcing} partial order $\F_\kappa$ are partial functions
$p\from\kappa\to\kappa$ such that $|p|<\kappa$ and if $\gamma\in\dom(p)$, then $\gamma$ is inaccessible,
$p\image\gamma\of\gamma$ and $|p\restrict\gamma|<\gamma$. We order the conditions by reverse inclusion. If
$p$ is a condition in $\F_\kappa$ for which $p(\gamma)=\alpha$, then $p$ naturally splits into two pieces,
$p=p_0\union p_1$, where $p_0=p\restrict\gamma$ and $p_1=p\restrict[\gamma,\kappa)$. Since
$p\image\gamma\of\gamma$, it follows that $p_0$ is in the fast function forcing partial order $\F_\gamma$
at $\gamma$. By splitting every stronger condition $q$ at $\gamma$ in the same way, namely $q\mapsto
q_0\union q_1$, one concludes that the forcing $\F_\kappa\restrict p$ below the condition $p$ is isomorphic
to $(F_\gamma\restrict p_0)\cross (\F_{\gamma,\kappa}\restrict p_1)$, where $\F_{\gamma,\kappa}$ is the set
of conditions in $\F_\kappa$ with domain in the interval $[\gamma,\kappa)$. By taking the unions of
conditions, one can see that $\F_{\gamma,\kappa}$ is $\leqgamma$-directed closed. In particular, if
$p=\{\<\gamma,\alpha>\}$ is the condition consisting of a function with one element in its domain, then
$\F_\kappa\restrict p\iso \F_\gamma\cross\F_{\gamma,\kappa}\restrict p$.

Let me first prove the basic theorem in the easiest case, that of a measurable cardinal. After this, I will
explain how to adapt the proof to the other large cardinal notions.

\Theorem. If\/ $\kappa$ is measurable, then there is a forcing extension preserving this in which
$\LDmeas_\kappa$ holds.\label{FFFmeas}

\Proof: This is essentially \cite[Theorem 2.3 ]{Hamkins2000:LotteryPreparation}. Suppose that $\kappa$ is
measurable in $V$. Without loss of generality, by forcing if necessary, I may assume that
$2^\kappa=\kappa^\plus$, since this can be forced without adding subsets to $\kappa$, and therefore without
destroying the measurability of $\kappa$. Let $f\from\kappa\to\kappa$ be $V$-generic for the fast function
forcing poset $\F_\kappa$. I claim that $\kappa$ is measurable in $V[f]$ and $f$ is an ordinal anticipating
Laver function there. To see this, fix any normal ultrapower $j:V\to M$ and any $\alpha<j(\kappa)$.
Consider the condition $p=\{\<\kappa,\alpha>\}$ in the poset $j(\F_\kappa)$. Below this condition, the
forcing $j(\F_\kappa)$ factors as $\F_\kappa\cross\F_{\kappa,j(\kappa)}\restrict p$. Observe that
$F_{\kappa,j(\kappa)}$ has size $j(\kappa)$ in $M$, and therefore has at most $j(2^\kappa)$ many open dense
sets in $M$. But $|j(2^\kappa)|^V\leq|\{h:\kappa\to 2^\kappa\st h\in V\}|^V=(2^\kappa)^\kappa=\kappa^\plus$
in $V$, so we may enumerate these dense sets in $V$ in a $\kappa^\plus$ sequence
$\<D_\alpha\st\alpha<\kappa^\plus>$. And since $M^\kappa\of M$, every initial segment of this enumeration
is in $M$. By diagonalizing this sequence, we may construct in $V$ a descending sequence of conditions
$\<p_\alpha\st\alpha<\kappa^\plus>$ below $p$ such that $p_\alpha\in D_\alpha$. At successor stages, simply
meet the next dense set; at limit stages, the union of the previous conditions is in $M$ because
$M^\kappa\of M$, and then one can meet the dense set below this condition. Thus, the function
$\ftail=\union_\alpha p_\alpha$ is in $V$, but is $M$-generic for the forcing $\F_{\kappa,j(\kappa)}$ below
$p$. Thus, by reversing the isomorphism, the function $f\union \ftail$ is $M$-generic for the forcing
$j(\F_\kappa)$. Since $j$ fixes every element of $\F_\kappa$, it follows that $j$ lifts to $j:V[f]\to
M[j(f)]$, where $j(f)=f\union\ftail$ and $j(\tau_f)=j(\tau)_{j(f)}$ for any $\F_\kappa$-name $\tau$. Since
$\kappa$ remains the critical point of this lifted embedding, we may conclude that $\kappa$ is measurable
in $V[f]$; and since $j(f)(\kappa)=\ftail(\kappa)=p(\kappa)=\alpha$, we see that $f$ anticipated $\alpha$.

It remains to check that there is an actual set-anticipating Laver function in $V[f]$. Fix any enumeration
of $V_\kappa$ in order type $\kappa$: $\vec a=\<a_\alpha\st\alpha<\kappa>$. Define
$\ell(\gamma)=(a_{f(\gamma)})_{f\restrict\gamma}$, provided that this makes sense, i.e., provided that
$\gamma$ is in the domain of $f$, that $f\restrict\gamma$ is $V$-generic for $\F_\gamma$, and that
$a_{f(\gamma)}$ is an $\F_\gamma$-name. I claim that $\ell$ is a $\LDmeas_\kappa$ Laver function. To see
this, fix any set $A\in H(\kappa^\plus)$. This set has an $\F_\kappa$-name $\dot A$ in $V$, and we may
assume $\dot A\in H(\kappa^\plus)^V$. Fix any normal ultrapower embedding $j:V\to M$ in $V$. Since
$M^\kappa\of M$ in $V$, it follows that $\dot A\in M$, and consequently $\dot A$ appears in the enumeration
$j(\vec a)$ of $M_{j(\kappa)}$ in $M$. So $\dot A=j(\vec a)(\alpha)$ for some particular $\alpha$. By the
lifting argument given previously, we may lift the embedding $j:V\to M$ to an embedding $j:V[f]\to M[j(f)]$
for which $j(f)(\kappa)=\alpha$. For this embedding, therefore, $j(\ell)(\kappa)$ is obtained by taking the
$\alpha^\th$ element in the sequence $j(\vec a)$, that is, $\dot A$, and interpreting it with the generic
object $j(f)\restrict\kappa=f$. Thus, $j(\ell)(\kappa)=(\dot A)_f=A$, as desired.\QED

A modified version of this forcing will guarantee that the fast function $f$ has a uniform nature, so that
for every measurable cardinal $\gamma<\kappa$, the restriction $f\restrict\gamma$ is generic for the
forcing up to $\gamma$. More generally, I will propose a modified class version of fast function forcing
that will produce a universal $\LDmeas$ Laver function that works with every measurable cardinal. In
addition, by applying Theorem \ref{Fast Function Flexibility Theorem} to embeddings arising from the
iteration of a normal measure, one can see that the fast function $f$ is able to anticipate any ordinal at
all; that is, for any ordinal $\alpha$ there is an embedding $j:V[f]\to M[j(f)]$ with critical point
$\kappa$ such that $j(f)(\kappa)=\alpha$.

\Theorem. There is a (class) forcing extension, preserving all measurable cardinals (and creating no new
measurable cardinals) in which there is a universal $\LDmeas$ Laver function.\label{FFFmeasUniform}

\Proof: The basic strategy is first to add a class fast function, whose restriction to any measurable
cardinal will be an ordinal anticipating Laver function, and then from this function to build as before a
set-anticipating Laver function. We may assume without loss of generality, by forcing if necessary, that
the \GCH\ holds, since it is known how to force this while preserving all measurable cardinals (and the
results on gap forcing in \cite{Hamkins2001:GapForcing} guarantee that this forcing creates no new
measurable cardinals). By further class forcing, if necessary, we may assume \AC\ for classes (this forcing
adds no new sets at all). Consider now a modified class forcing version $\F^*$ of fast function forcing:
conditions are partial (set) functions $p\from\ORD\to\ORD$ such that (1) every $\gamma\in\dom(p)$ is
inaccessible, with $p\image\gamma\of\gamma$ and $|p\restrict\gamma|<\gamma$, (2) if $\delta$ is
inaccessible and $p\image\delta\of\delta$, then $|p\restrict\delta|<\delta$, and (3) if $\kappa$ is any
measurable cardinal, then $p\image\kappa\of\kappa$ and $|p\restrict\kappa|<\kappa$. The last modification,
in effect, slows the functions down; they cannot jump over the next measurable cardinal. The conditions are
ordered by reverse inclusion. It is not difficult to see that for any measurable cardinal $\kappa$, the
forcing $\F^*$ factors as $\F^*_\kappa\cross\F^*_{\kappa,\infty}$, where $\F^*_\kappa$ is the set of
conditions in $\F^*$ with domain (and range) in $\kappa$, and $\F^*_{\kappa,\infty}$ is those with domain
in $[\kappa,\infty)$. With the modified forcing $\F^*$, there is no need to factor only below a condition,
since no condition can hop over $\kappa$.

Suppose now that $f^*\from\ORD\to\ORD$ is $V$-generic for the forcing $\F^*$. I claim that $f^*$ is a
universal ordinal anticipating Laver function. To see this, suppose that $\kappa$ is a measurable cardinal
and consider the restriction $f=f^*\restrict\kappa$. It is easy to see that $f$ is $V$-generic for the
forcing $\F^*_\kappa$, and the remaining forcing $F^*_{\kappa,\infty}$ is highly closed, definitely adding
no new $\kappa$ sequences over $V[f]$. In particular, if
$\ftail^*=f^*_{\kappa,\infty}=f^*\restrict[\kappa,\infty)$, then $\kappa$ remains measurable in
$V[\ftail^*]$. I will now argue that $f$ is an ordinal anticipating Laver function in $V[f^*]$. For this,
we proceed just as in Theorem \ref{FFFmeas}: fix any ultrapower embedding $j:V[\ftail^*]\to M[j(\ftail^*)]$
by a normal measure in $V$, and consider the forcing $j(\F^*_\kappa)$, which factors as
$\F^*_\kappa\cross\F^*_{\kappa,j(\kappa)}$, below the condition $p=\{\<\kappa,\alpha>\}$. Below this
condition in $\F^*_{\kappa,j(\kappa)}$, build in $V$ by diagonalization an $M$-generic filter $\ftail$, and
lift the embedding to $j:V[\ftail^*][f]\to M[j(\ftail^*)][j(f)]$, where $j(f)=f\union\ftail$. Since
$\ftail$ included $p$, it follows that $j(f)(\kappa)=\ftail(\kappa)=p(\kappa)=\alpha$, as desired.

Next, I produce from $f^*$ a universal set anticipating Laver function. Using \AC\ for classes, fix any
class well-ordering of the universe $V$ in order type $\ORD$. It is a simple exercise to produce from this
an order $\trianglelt$ for which the objects of $V_\kappa$ appear as the first $\kappa$ many elements of
the order, whenever $\kappa$ is an inaccessible cardinal. Next, as in Theorem \ref{FFFmeas}, we construct
the Laver function by defining $\ell(\gamma)$ to be $(a_{f(\gamma)})_{f^*\restrict\gamma}$, where $a_\beta$
is the $\beta^\th$ object in the well order $\trianglelt$, provided that $a_{f(\gamma)}$ is a
$\F^*_\gamma$-name. I claim that $\ell\from\ORD\to V[f^*]$ is a universal Laver function. To see this,
suppose that $\kappa$ is any measurable cardinal and $A\in H(\kappa^\plus)^{V[f^*]}$. Since the forcing
above $\kappa$ is highly closed, we know that actually $A\in H(\kappa^\plus)^{V[f^*\restrict\kappa]}$, and
so $A$ has a $\F^*_\kappa$-name $\dot A$ in $H(\kappa^\plus)^V$, that is, $A=\dot A_f$. Fix any normal
ultrapower $j:V\to M$ in $V$. As in Theorem \ref{FFFmeas}, the name $\dot A$ is the $\beta^\th$ object of
the well order $j(\trianglelt)$ for some $\beta<j(\kappa)$, and the argument of the previous paragraph
explains how to lift this embedding to $j:V[f^*]\to M[j(f^*)]$ such that $j(f)(\kappa)=\beta$, that is,
$j(f^*)(\kappa)=\beta$. From this, it follows that $j(\ell)(\kappa)=(\dot A)_{j(f^*)\restrict\kappa}=\dot
A_f=A$, as desired. So $\ell$ is a universal $\LDmeas$ Laver function.

Finally, I complete the argument by pointing out the Gap Forcing Theorem of \cite{Hamkins2001:GapForcing}
shows that this forcing creates no new measurable cardinals, since the GCH forcing followed by $\F^*$
admits a very low gap.\QED

These proofs adapt easily to the case of a weakly compact cardinal.

\Theorem. If $\kappa$ is a weakly compact cardinal, then there is a forcing extension in which
$\LDwc_\kappa$ holds. Indeed, there is a class forcing extension, preserving all weakly compact cardinals,
in which there is a universal $\LDwc$ Laver function.\label{FFFwc}

\Proof: The first part of this theorem is essentially \cite[Theorem 2.4]{Hamkins2000:LotteryPreparation}.
To handle only one weakly compact cardinal, we simply perform Woodin's fast function forcing $\F_\kappa$ as
in Theorem \ref{FFFmeas}. In particular, there is no need for the \GCH\ in this argument, since it is
considerably easier to lift the small weakly compact embeddings than full class embeddings on $V$. Suppose
that $f\from\kappa\to\kappa$ is a fast function. I claim that for any nice structure $M$ in $V$ and any
canonical weakly compact embedding $j:M\to N$, so that $N$ is a transitive structure of size $\kappa$ and
$N^\ltkappa\of N$ in $V$, there is a lift of the embedding to $j:M[f]\to N[j(f)]$. This is because the
forcing $\F_\kappa$ is definable from $\kappa$ in $M$, and the forcing $j(\F_\kappa)$ factors below the
condition $p=\{\<\kappa,\alpha>\}$ as $\F_\kappa\cross\F_{\kappa,j(\kappa)}$, for any $\alpha<j(\kappa)$.
Since there are only $\kappa$ many dense sets in $N$ (all of $N$ has size $\kappa$), and the forcing
$\F_{\kappa,j(\kappa)}$ is $\leqkappa$-closed, the diagonalization argument of Theorem \ref{FFFmeas} easily
produces an $N$-generic filter $\ftail\of\F_{\kappa,j(\kappa)}$ below $p$. And since $f$ is $V$-generic for
$\F_\kappa$, it is surely also $N[\ftail]$-generic, so $f\union\ftail$ is $N$-generic for $j(\F_\kappa)$
below $p$. The embedding therefore lifts to $j:M[f]\to N[j(f)]$ where $j(f)=f\union\ftail$. In particular,
$j(f)(\kappa)=\ftail(\kappa)=p(\kappa)=\alpha$. This almost shows that $f$ is an ordinal anticipating
weakly compact Laver function, except that so far we have only considered structures of the form $M[f]$,
where $M$ is a nice structure in $V$. If $M'$ is any nice structure of size $\kappa$ in $V[f]$, then $M'$
has an $\F_\kappa$ name $\dot M'$ in $V$ of size $\kappa$. And this name $\dot M'$ can be placed inside a
nice structure $M$ of size $\kappa$ in $V$. Thus, since $\dot M'$ and $f$ are both in $M[f]$, it follows
that $M'=\dot M'_f$ is in $M[f]$. If $j:M[f]\to N[j(f)]$ was the embedding I constructed above, then the
restriction $j\restrict M':M'\to j(M')$ provides a witnessing embedding for $M'$. In particular, if $f\in
M'$, then this restriction also has $j(f)(\kappa)=\alpha$, so $f$ is an ordinal anticipating weakly compact
Laver function in $V[f]$. One now constructs the set anticipating Laver function $\ell$ just as before: fix
a well ordering $\trianglelt$ of $V_\kappa$ such that $V_\gamma$ is enumerated in order type $\gamma$ for
any inaccessible cardinal $\gamma$, and let $\ell(\gamma)$ be $(a_{f(\gamma)})_{f\restrict\gamma}$, where
$a_\beta$ is the $\beta^\th$ element in the order, provided that $a_{f(\gamma)}$ is indeed a $\F_\gamma$
name. If $A\in H(\kappa^\plus)^{V[f]}$, then $A=\dot A_f$ for some name $\dot A\in H(\kappa^\plus)^V$.
Suppose that $M$ is a nice structure of size $\kappa$ in $V$ such that $f$, $\dot A$ and $\trianglelt$ are
in $M$, and fix any weakly compact embedding $j:M\to N$ in $V$. The object $\dot A$ is in $N$, and
therefore is the $\beta^\th$ element in the order there. The lifting argument above allows us to lift the
embedding to $j:M[f]\to N[j(f)]$ such that $j(f)(\kappa)=\beta$. It follows that
$j(\ell\restrict\kappa)(\kappa)=\dot A_{j(f)\restrict\kappa}=\dot A_f=A$, as required. So $\ell$ is a
$\LDwc_\kappa$ Laver function in $V[f]$.

The universal result is obtained as in Theorem \ref{FFFmeasUniform} by slowing down the fast function
forcing so as to respect every weakly compact cardinal. Specifically, in the weakly compact context we must
ensure that every weakly compact cardinal has a chance to be in the domain of the function. Let $\F^*$
consist of the conditions $p\from\ORD\to\ORD$ such that, as before, (1) every $\gamma\in\dom(p)$ is an
inaccessible cardinal with $p\image\gamma\of\gamma$ and $|p\restrict\gamma|<\gamma$, (2) if $\delta$ is
inaccessible and $p\image\delta\of\delta$, then $|p\restrict\delta|<\delta$, and additionally (3) whenever
$\gamma$ is a weakly compact cardinal, then $p\image\gamma\of\gamma$ and $|p\restrict\gamma|<\gamma$. Every
weakly compact cardinal will be closed under the resulting generic function $f^*\from\ORD\to\ORD$, and it
is easy to show as above that $f^*$ is a universal ordinal anticipating weakly compact Laver function. The
corresponding function $\ell$ will be a universal $\LDwc$ function in $V[f^*]$.\QED

\Theorem. If\/ $\kappa$ is an unfoldable cardinal, then there is a forcing extension in which
$\LDunf_\kappa$ holds. Indeed, there is a class forcing extension, preserving all unfoldable cardinals, in
which there is a universal $\LDunf$ Laver function.\label{FFFunf}

\Proof: This argument follows the standard lifting technique for unfoldability embeddings, given in
\cite{Hamkins2001:UnfoldableCardinals}, an adaptation of Woodin's factor technique for strongness
embeddings to the much smaller domains.

Suppose that $\kappa$ is unfoldable and $f\from\kappa\to\kappa$ is $V$-generic for Woodin's fast function
forcing $\F_\kappa$.  Fix any $\theta$ and any ground model canonical $(\theta+1)$-unfoldability embedding
$j:M\to N$, where $M$ is a nice structure of size $\kappa$ in $V$. That is, $j(\kappa)>\theta$,
$N^\ltkappa\of N$ and $N=\{\,j(h)(\beta_0,\beta_1,\ldots,\beta_n)\st h\in
M\And\kappa\leq\beta_0\leq\cdots\leq\beta_n\leq\theta\,\}$. Fix any $\alpha<j(\kappa)$, and let
$X=\set{j(h)(\kappa,\alpha,\theta)\st h\in M}$ be the hull of the pair $\<\kappa,\theta>$ in $N$. By
verifying the Tarski-Vaught criterion, one deduces that $X\elesub N$. Below the condition
$p=\{\<\kappa,\alpha>,\<\kappa',\theta>\}$, where $\kappa'$ is the next inaccessible in $N$ above $\kappa$
and $\alpha$, the forcing $j(\F_\kappa)$ factors as $\F_\kappa\cross(\F_{\kappa,j(\kappa)}\restrict p)$,
and $\F_{\kappa,j(\kappa)}\restrict p$ is $\leqtheta$-closed in $N$, and hence also in $X$. Since $X$ has
size $\kappa$, it has only $\kappa$ many dense sets for $\F_{\kappa,j(\kappa)}$, and so we may by
diagonalization construct an $X$-generic filter $\ftail^0\of\F_{\kappa,j(\kappa)}\intersect X$ by simply
lining up the dense sets in a $\kappa$-sequence and meeting them one-by-one. I claim that the filter
$\ftail$ generated by $\ftail^0$ in $\F_{\kappa,j(\kappa)}$ is fully $N$-generic. To see this, suppose that
$D\in N$ is an open dense subset of $\F_{\kappa,j(\kappa)}$. By our assumption on the embedding $j$, we
know that $D=j(\vec D)(\beta_0,\ldots,\beta_n)$ for some function $\vec D\in M$ and ordinals
$\kappa\leq\beta_0\leq\cdots\leq\beta_n\leq\theta$. Let $\bar D$ be the intersection of all $j(\vec
D)(\alpha_0,\ldots,\alpha_n)$, over $\kappa\leq\alpha_0\leq\cdots\leq\alpha_n\leq\theta$, such that this is
an open dense subset of $\F_{\kappa,j(\kappa)}$. Since $\bar D$ is definable from $j(\vec D)$, $\kappa$ and
$\theta$, it follows that $\bar D$ is in $X$. Furthermore, since $\F_{\kappa,j(\kappa)}$ is
$\leqtheta$-closed below $p$, and there are at most $\theta$ many open dense sets in this intersection, it
follows that $\bar D$ remains open and dense in $F_{\kappa,j(\kappa)}$ below $p$. Thus, since it is in $X$,
the filter $\ftail^0$ meets $\bar D$, and consequently $\ftail$ meets it as well. Finally, since $\bar D\of
D$, we conclude that $\ftail$ meets $D$, as desired. So $\ftail$ is $N$-generic, and we may lift the
embedding to $j:M[f]\to N[j(f)]$, where $j(f)=f\union\ftail$. Since
$j(f)(\kappa)=\ftail(\kappa)=p(\kappa)=\alpha$, we know that $f$ is an ordinal anticipating Laver function
in $V[f]$. From this function, one builds the $\LDunf_\kappa$ Laver function $\ell$ as in Theorem
\ref{FFFwc}.

The universal version of this theorem is proved just as in Theorem \ref{FFFwc}. Here, we modify the fast
function forcing to include only those conditions $p$ with $|p\restrict\kappa|<\kappa$ and
$p\image\delta\of\delta$ for any unfoldable cardinal $\delta$, thereby producing the modified forcing
$\F^*$. If $f^*$ is $V$-generic for $\F^*$, then $f=f^*\restrict\kappa$ is $V$-generic for $\F^*_\kappa$
for any unfoldable cardinal $\kappa$. One then argues that this is an ordinal anticipating unfoldability
Laver function at $\kappa$ in $V[f]$. From this function, one builds the set anticipating $\LDunf$ Laver
function $\ell$ as before. The remainder of the forcing $\F^*_{\kappa,\infty}$ is $\leq$-closed, and
therefore creates no new transitive structures of size $\kappa$, so $f$ and $\ell$ retain these properties
in $V[f^*]$, and the proof is complete.\QED

A similar argument works in the case of $\Pi^m_n$-indescribability embeddings, to obtain the following.

\Theorem. If $\kappa$ is $\Pi^m_n$-indescribable, then there is a forcing extension where
$\LD^{\hbox{\!\!\tiny$\Pi^m_n$-ind}}_\kappa$ holds. Indeed, there is a class forcing extension, preserving
all $\Pi^m_n$-indescribable cardinals, in which there is a universal $\LD^{\hbox{\!\!\tiny$\Pi^m_n$-ind}}$
Laver function.

And the construction works for strongly unfoldable cardinals as well.

\Theorem. If\/ $\kappa$ is strongly unfoldable, then there is a forcing extension in which $\LDsunf_\kappa$
holds.

\Proof: I will adapt the unfoldability argument of Theorem \ref{FFFunf} to the case of strong
unfoldability. Suppose in that argument that $\kappa$ is strongly unfoldable, and that we have added the
fast function $f\from\kappa\to\kappa$. Suppose in that argument that the ground model embedding $j:M\to N$
is not merely a $\theta$-unfoldability embedding, but a $\theta$-strong unfoldability embedding, so that we
know in addition that $V_\theta\of N$ and $N=\{\,j(h)(\beta_0,\ldots,\beta_n)\st h\in
M\And\kappa\leq\beta_0\leq\cdots\leq\beta_n\leq\beth_\theta\,\}$ (one needs $\beth_\theta$ here in order to
capture all of $V_\theta$). The argument of Theorem \ref{FFFunf} shows how to lift the embedding to
$j:V[f]\to M[j(f)]$, such that $j(f)(\kappa)=\alpha$, where $\alpha$ is any desired ordinal below
$j(\kappa)$. In that argument, the lifted embedding is a $\theta$-unfoldability embedding. I claim that
since we began here with a $\theta$-strong unfoldability embedding, this lift is actually a $\theta$-strong
unfoldability embedding. This is because from $V_\theta\of N$ it follows that $(V[f])_\theta\of N[f]\of
N[j(f)]$. And since $M[f]$ and $M'$ have the same $V_\kappa$, it follows that $N[j(f)]$ and $j(M')$ have
agree up to rank $j(\theta)$, so $(V[f])_\theta\of j(M')$ as well. Thus, the restriction $j\restrict
M':M'\to j(M')$ is also a $\theta$-strong unfoldability embedding on $M'$, for which $j(f)(\kappa)=\alpha$.
Now, from the ordinal anticipating Laver function $f$ one can easily build a set anticipating Laver
function as before, so the proof is complete.\QED

One naturally wonders whether the universal form $\LDsunf$ is forceable, while preserving all strongly
unfoldable cardinals. By following the previous argument, using the modified universal fast function
forcing $\F^*$ consisting of conditions $p\in\F$ that satisfy $|p\restrict\kappa|<\kappa$ and
$p\image\kappa\of\kappa$ for any strongly unfoldable cardinal $\kappa$, one can see that the universal
function $f^*$ obtained will almost be an ordinal anticipating strong unfoldability Laver function. For any
strongly unfoldable cardinal $\kappa$, any ordinal $\theta$ and any $\alpha$ up to the next inaccessible
cluster point of $\dom(f^*)$ above $\kappa$, the ground model canonical $\theta$-strong unfoldability
embeddings will lift to $j:M[f]\to N[j(f)]$ with $j(f)(\kappa)=\alpha$, such that $(V[f^*])_\theta\of
N[j(f)]$. The problem arises when $\alpha$ is too large, so that in the forcing in $N$ between $\kappa$ and
$\theta$ we cannot use $f^*\restrict(\kappa,\theta)$, because the master condition $p=\{\<\kappa,\alpha>\}$
jumps too high. Thus, although we will be able to lift the embedding, it will not include enough of $f^*$
to be a $\theta$-strong unfoldability embedding. What is needed is to simultaneously force some
indestructibility for strong unfoldability, as in the argument of Theorem \ref{ForcingLDsuperstrong}, so
that this generic object can be recovered. Such kind of indestructibility arguments for strong
unfoldability will be the focus of a forthcoming paper.

\Theorem. If $\kappa$ is a Ramsey cardinal, then there is a forcing extension preserving this in which
$\LDram_\kappa$ holds. Indeed, there is a class forcing extension, preserving all Ramsey cardinals, in
which there is a universal $\LDram$ Laver function.

\Proof: Let me first handle just one Ramsey cardinal $\kappa$. Let $f\from\kappa\to\kappa$ be a $V$-generic
fast function, and fix a well ordering $\trianglelt$ of $V_\kappa$ in order type $\kappa$ in $V$. Let
$\ell(\gamma)$ be $\tau_{f\restrict\gamma}$, where $\tau$ is the $f(\gamma)^\th$ element of $V_\kappa$ with
respect to $\trianglelt$, provided that this object is in fact an $\F_\gamma$-name. I claim that $\ell$ is
a $\LDram_\kappa$ Laver function in $V[f]$. Fix any nice structure $M$ of size $\kappa$ in $V[f]$ with
$\ell\in M$ and any $A\in H(\kappa^\plus)^M$. The set $A$ has a name $\dot A$ of size $\kappa$ in $V$. The
structure $M$ also has a name $\dot M$ of size $\kappa$ in $V$, and we may find a nice structure $M_0$ of
size $\kappa$ in $V$ with $\dot A, \dot M$ and $\trianglelt$ in $M_0$. Since $\kappa$ is a Ramsey cardinal
in $V$, there is a weakly amenable $M_0$-normal measure $F$ in $V$ such that $\<M_0,F>$ is iterable. Let
$j:M_0\to M_1$ be the ultrapower of $M_0$ by $F$. It is easy to see that $\dot A$ is in $M_1$, and so it is
the $\beta^\th$ object with respect to $j(\trianglelt)$. Below the condition $\<\kappa,\beta>$, the forcing
$j(\F)$ factors as $\F\cross\Ftail$, where $\Ftail$ is $\leqkappa$-closed in $M_1$. Since $M_1^\ltkappa\of
M_1$ and there are at most $\kappa$ many dense subsets of $\Ftail$ in $M_1$, we may diagonalize to
construct in $V$ an $M_1$ generic object $\ftail$, and lift the embedding to $j:M_0[f]\to M_1[j(f)]$, with
$j(f)(\kappa)=\beta$. It follows that $j(\ell)(\kappa)=\dot A_f=A$, so the function $\ell$ anticipates $A$,
as desired. The lifted embedding remains iterable, because the filter induced by the lift is
$\ltkappa$-closed in $V[f]$, as $M[f]$ is $\ltkappa$-closed in $V[f]$.

One obtains the universal form as before, by performing a modified class forcing $\F$, consisting of fast
function forcing conditions $p$ such that if $\kappa$ is a Ramsey cardinal, then $p\image\kappa\of
V_\kappa$ and $\card{p\restrict\kappa}<\kappa$. For any particular Ramsey cardinal $\kappa$, this forcing
factors as $\F_\kappa\cross\F_{\kappa,\infty}$. The cardinal $\kappa$ remains Ramsey in
$V^{\F_{\kappa,\infty}}$, because this is $\leqkappa$-directed closed forcing, and the argument of the
previous paragraph shows that subsequent forcing with $\F_\kappa$ produces the $\LDram_\kappa$ Laver
function. Finally, the arguments of \cite{Hamkins2001:GapForcing} show that no new Ramsey cardinals are
created.\QED

Let me now turn to the larger large cardinals. The fact is that fast function forcing produces Laver
functions in a very general way. To support of this view, I will make use of the following two results from
previous work.

\Theorem Fast Function Flexibility Theorem.{(\cite[Theorem 1.11]{Hamkins2000:LotteryPreparation})} Suppose
that $f\from\kappa\to\kappa$ is a fast function added generically over $V$ and that $j:V[f]\to M[j(f)]$ is
an embedding (either internal or external to $V[f]$) with critical point $\kappa$. Then for any
$\alpha<j(\kappa)$ there is another embedding $j^*:V[f]\to M[j^*(f)]$ such that:
 \begin{enumerate}
 \item $j^*(f)(\kappa)=\alpha$,
 \item $j^*\restrict V=j\restrict V$,
 \item $M[j^*(f)]\of M[j(f)]$, and
 \item If $\alpha$ is less than the next inaccessible cluster point of $\dom(j(f))$ beyond $\kappa$,
        then $M[j^*(f)]=M[j(f)]$. In
            this case, if $j$ is the ultrapower by a standard measure $\eta$ concentrating on a set in $V$,
            then $j^*$ is the ultrapower by a standard measure $\eta^*$ concentrating on the
            same set, and moreover $\eta\intersect V=\eta^*\intersect V$ and $[\id]_\eta=[\id]_{\eta^*}$.
 \end{enumerate}

A measure is {\df standard} if the critical point $\kappa$ of the corresponding ultrapower embedding
$j:V[f]\to M[j(f)]$ is definable in $M[j(f)]$ from $s=[\id]_\mu$ and parameters in $\ran(j\restrict V)$.
Thus, any normal measure on $\kappa$ is standard (since in this case $s=\kappa$), as is any
supercompactness measure (since in this case $\kappa$ is the least element not in $s=j\image\theta$). Also,
\cite[Lemma 2.7] {Hamkins2000:LotteryPreparation} shows that in the types of forcing extensions of this
paper, every strong compactness measure is isomorphic to a standard strong compactness measure.

\Theorem Generalized Laver Function Theorem.{(\cite[Theorem 2.2]{Hamkins2000:LotteryPreparation})} After
fast function forcing $V[f]$, there is a function $\ell\from\kappa\to (V[f])_\kappa$ with the property that
for any embedding $j:V[f]\to M[j(f)]$ with critical point $\kappa$ (whether internal or external) and for
any $z\in H(\lambda^\plus)^{M[j(f)]}$, where $\lambda=j(f)(\kappa)$, there is another embedding
$j^*:V[f]\to M[j(f)]$ such that:
 \begin{enumerate}
  \item $j^*(\ell)(\kappa)=z$,
  \item $M[j^*(f)]=M[j(f)]$,
  \item $j^*\restrict V=j\restrict V$, and
  \item If $j$ is the ultrapower by a standard measure $\eta$ concentrating on a set in $V$, then $j^*$
  is the ultrapower by a standard measure $\eta^*$ concentrating on the same set and moreover
  $\eta^*\intersect V=\eta\intersect V$ and  $[\id]_{\eta^*}=[\id]_\eta$.
 \end{enumerate}

The function $\ell$ is defined from $f$ in the same way as in Theorem \ref{FFFmeas}. The point now is that
this general theory of fast function forcing allows one to deduce $\LDstar_\kappa$ in a great variety of
cases. All one really needs to know is that fast function forcing preserves the large cardinal in question,
so that there are embeddings $j:V\to M$ as in Theorem \ref{Fast Function Flexibility Theorem}, and then
Theorem \ref{Generalized Laver Function Theorem} essentially says that $\LDstar_\kappa$ holds, provided
that one can deduce that $j^*$ is the appropriate type of embedding (or has an appropriate factor).

\Corollary. If $\kappa$ is strongly compact, there is a forcing extension in which $\LDstrc_\kappa$ holds.

\Proof: Theorem 1.7 of \cite{Hamkins2000:LotteryPreparation} show that fast function forcing preserves
every strongly compact cardinal $\kappa$. Let $\ell$ be the function of Theorem \ref{Generalized Laver
Function Theorem}. By \cite[Theorem 1.12]{Hamkins2000:LotteryPreparation}, there is for any $\theta$ a
$\theta$-strong compactness embedding $j:V[f]\to M[j(f)]$ for which $j(f)(\kappa)>\theta$. For any $z\in
H(\theta^\plus)^{V[f]}$, now, we may now apply Theorem \ref{Generalized Laver Function Theorem} to produce
an embedding $j^*:V[f]\to M[j^*(f)]$ for which $j^*(\ell)(\kappa)=z$. If $\theta=\theta^\ltkappa$, then
Theorem 2.6 of \cite{Hamkins2000:LotteryPreparation} shows that the embedding $j^*$ can be chosen to be a
$\theta$-strong compactness embedding. So $\ell$ witnesses $\LDstrc_\kappa$ in $V[f]$.\QED

\Corollary. If $\kappa$ is $2^{\theta^\ltkappa}$-strongly compact, then there is a forcing extension in
which $\LDthetastrc_\kappa$ holds.

\Proof: This is the amount of strong compactness that is used in the previous argument.\QED

\Corollary. If $\kappa$ is $\theta$-supercompact, then there is a forcing extension where
$\LDthetasc_\kappa$ holds.\label{ForcingLDscLocal}

\Proof: Suppose $\kappa$ is $\theta$-supercompact. Since this implies that $\kappa$ is also
$\theta^\ltkappa$-supercompact, we may replace $\theta$ with $\theta^\ltkappa$ if necessary and assume
$\theta^\ltkappa=\theta$. Now, by forcing if necessary, we may assume $2^\theta=\theta^\plus$, since the
forcing to accomplish this adds no subsets of $P_\kappa\theta$. Next, Theorem 1.10 of
\cite{Hamkins2000:LotteryPreparation} shows that fast function forcing preserves the $\theta$
supercompactness of $\kappa$, so there are embeddings $j:V[f]\to M[j(f)]$ as in Theorem \ref{Fast Function
Flexibility Theorem}. Necessarily,  by the proof of \cite[Theorem 1.13]{Hamkins2000:LotteryPreparation},
the next inaccessible cluster point of $\dom(j(f))$ beyond $\kappa$ is at least $\theta$. Thus, by Theorem
\ref{Fast Function Flexibility Theorem} we may find an embedding for which $j(f)(\kappa)>\theta$. Finally,
Theorem \ref{Generalized Laver Function Theorem} implies that the function $\ell$ witnesses
$\LDthetasc_\kappa$ in $V[f]$, as desired.\QED

This theorem shows that the consistency strength of $\LDthetasc$ is not greater than $\kappa$ being
$\theta$-supercompact. Consequently, the hypothesis of Corollary \ref{LaverLocal} is not optimal in terms
of consistency strength. It is open whether every $\theta$-supercompact cardinal $\kappa$ has a
$\LDthetasc_\kappa$ Laver function, at least when $\kappa<\theta$. When $\kappa=\theta$, the outright
implication is ruled out by Theorem \ref{LDmeasFailsInL[mu]}.

\Corollary.If $\kappa$ is a $\theta$-strong cardinal and $\theta$ is either a successor ordinal or has
cofinality at least $\kappa$, then there is a forcing extension in which $\LDthetastr$
holds.\label{ForcingLDstrLocal}

\Proof: The standard reverse Easton iteration arguments establish that one may force
$2^\kappa=\kappa^\plus$ while preserving the $\theta$-strongness of $\kappa$. After this, \cite[Theorem
1.6]{Hamkins2000:LotteryPreparation} shows that fast function forcing preserves the $\theta$-strongness of
$\kappa$ and provides a ordinal-anticipating Laver function. A full set-anticipating $\LDthetastr_\kappa$
Laver function can be constructed as in Theorem \ref{ForcingLD} above.\QED

One does not expect fast function forcing up to the cardinal in question by itself to produce Laver
functions for superstrong, extendible, almost huge or huge embeddings, since these embeddings have
$V_{j(\kappa)}\of M$, and consequently the function $j(f)$ would have to include fully $V$-generic portions
above $\kappa$, but there are no such objects in $V[f]$.  Rather, for these kind of embeddings, one expects
to force above $\kappa$. This introduces problems when one wants the fast function to anticipate an ordinal
larger than the value selected above $\kappa$ in $V$, and to overcome this difficulty one must also force a
kind of indestructibility that allows for this intermediate forcing to be recovered. As a sample of this
kind of argument, let me consider the case of superstrong cardinals.

\Theorem. If $\kappa$ is superstrong, then this is preserved to a forcing extension in which
$\LD_\kappa^{\hbox{\!\!\tiny superstrong}}$ holds.\label{ForcingLDsuperstrong}

\Proof: Suppose that $\kappa$ is superstrong, with embedding $j:V\to M$, so that $V_{j(\kappa)}\of M$ and
$M=\set{j(h)(s)\st h:[\kappa]^\ltomega\to V\And s\in [j(\kappa)]^\ltomega}$. Suppose that $f\from\kappa\to
\kappa$ is $V$-generic for fast function forcing $\F_\kappa$. Factor $j(\F_\kappa)$ as
$\F_\kappa\cross\F_{\kappa,j(\kappa)}$ below some condition $p=\{\<\kappa,\alpha>\}$, and let
$f_{\kappa,j(\kappa)}$ be $V[f]$-generic, so that $j$ lifts to $j:V[f]\to M[j(f)]$. Since
$\F_{\kappa,j(\kappa)}$ is $\leqkappa$-closed in $V$, it is $\leqkappa$-distributive in $V[f]$, and so the
extender embedding lifts further to $j:V[f][f_{\kappa,j(\kappa)}]\to M[j(f)][j(f_{\kappa,j(\kappa)})]$. Let
$\P_\kappa$ be the reverse Easton $\kappa$-iteration that at stage $\gamma\in\dom(f)$ forces with the
lottery sum of all posets $\Q\in H(f(\gamma)^\plus)$ that are $\leqgamma$-distributive. The forcing
$j(\P_\kappa)$ factors as $\P_\kappa*\P_{\kappa,j(\kappa)}$, where we opt for trivial forcing in the stage
$\kappa$ lottery. Let $G_\kappa*G_{\kappa,j(\kappa)}$ be $V[f][f_{\kappa,j(\kappa)}]$-generic for this
forcing. Once again, the embedding lifts fully to
$j:V[f][f_{\kappa,j(\kappa)}][G_\kappa][G_{\kappa,j(\kappa)}]\to
M[j(f)][j(f_{\kappa,j(\kappa)})][j(G_\kappa)][j(G_{\kappa,j(\kappa)})]$. Since $V_{j(\kappa)}\of M$ and all
the generic objects $f$, $f_{\kappa,j(\kappa)}$, $G_\kappa$ and $G_{\kappa,j(\kappa)}$ are present on the
$M$-side of the final embedding, this embedding witnesses that $\kappa$ remains superstrong in
$V[f][f_{\kappa,j(\kappa)}][G_\kappa][G_{\kappa,j(\kappa)}]$.

It remains to check that $\LD_\kappa^{\hbox{\!\!\tiny superstrong}}$ holds there. First let me check that
$f$ is an ordinal anticipating Laver function there. Certainly we may simply change the value of
$j(f)(\kappa)$ to be any value less than $\alpha$ without much change in the argument. The difficulty comes
when we want to hit some $\beta$ that is much bigger than $\alpha$. Suppose $\beta<j(\kappa)$. Let $\delta$
be the first element of $\dom f_{\kappa,j(\kappa)}$ with $\beta<\delta$. Reconsider the embedding $j:V\to
M$, but this time lift it to an embedding $j^*:V[f][f_{\kappa,j(\kappa)}]\to
M[j^*(f)][j^*(f_{\kappa,j(\kappa)})]$ where $j^*(f)(\kappa)=\beta$, so that $j^*(f)=f\union
f_{\delta,j(\kappa)}$. That is, since $j^*(f)(\kappa)$ jumps up to $\beta$, we simply omit the part of
$f_{\kappa,j(\kappa)}$ between $\kappa$ and $\delta$. Now, in the $j(\P_\kappa)$ forcing, instead of opting
for trivial forcing at stage $\kappa$, we reinsert the missing forcing by opting for the forcing
$\F_{\kappa,\delta}*\P_{\kappa,\delta}$, for which we have the appropriate generic objects. To summarize,
these maneuvers effectively reorganize the entire forcing $j(\F_\kappa)*j(\P_\kappa)$ as
$(\F_\kappa\cross\F_{\delta,j(\kappa)})*\P_\kappa*(\check
\F_{\kappa,\delta}*\P_{\kappa,\delta})*\P_{\delta,j(\kappa)}$. The embedding lifts as above to
$j^*:V[f][f_{\kappa,j(\kappa)}][G_\kappa][G_{\kappa,j(\kappa)}]\to
M[j^*(f)][j^*(f_{\kappa,j(\kappa)})][j^*(G_\kappa)][j^*(G_{\kappa,j(\kappa)})]$. Because the reorganized
forcing reinserts all the required generic objects on the $j$ side, this embedding is superstrong in
$V[f][f_{\kappa,j(\kappa)}][G_\kappa][G_{\kappa,j(\kappa)}]$. Furthermore, it has $j^*(f)(\kappa)=\beta$,
so $f$ is an ordinal anticipating superstrong Laver function. One may now construct a set anticipating
Laver function from $f$, as in Theorem \ref{FFFmeas}, and so $\LD_\kappa^{\hbox{\!\!\tiny superstrong}}$
holds in $V[f][f_{\kappa,j(\kappa)}][G_\kappa][G_{\kappa,j(\kappa)}]$, as desired.\QED

Let me conclude this section by showing a few alternative approaches to forcing the Laver diamond. The
first of these uses ordinary Silver forcing to add the set anticipating function directly, without need to
consider ordinal anticipating functions.

\Proof Alternative Proof of Theorem \ref{FFFmeas}: Suppose that $\kappa$ is a measurable cardinal and
$2^\kappa=\kappa^\plus$. Let $\P=\P_\kappa*\Q_\kappa$ be the ordinary Silver iteration of length
$\kappa+1$, that is, the $(\kappa+1)$-iteration which at every inaccessible stage $\gamma$ adds a Cohen set
to $\gamma$. This stage $\gamma$ forcing is isomorphic to the forcing which adds, by initial segments, a
function $f_\gamma:\gamma\to V[G_\gamma]_\gamma$. Suppose that $G*f$ is $V$-generic for
$\P_\kappa*\Q_\kappa$, where $f:\kappa\to V[G]_\kappa$. I claim that $f$ is a $\LDmeas_\kappa$ Laver
function in $V[G][f]$. To see this, suppose that $A\in H(\kappa^\plus)^{V[G][f]}$. The set $A$ must have a
name $\dot A\in V$ such that $A=\dot A_{G*f}$ and $\dot A\in H(\kappa^\plus)$. Fix any normal ultrapower
embedding $j:V\to M$. Consider the forcing $j(\P)$, which factors in $M$ as
$\P_\kappa*\Q_\kappa*\Ptail*j(\Q_\kappa)$, where $\Ptail$ is $\leqkappa$-closed in
$M^{\P_\kappa*\Q_\kappa}$. The filter corresponding to $G*f$ is $V$-generic and hence $M$-generic for the
first $\kappa+1$ many stages of this iteration. The forcing $\Ptail$ in $M[G][f]$ is $\leqkappa$-closed,
and so by diagonalization we may construct in $V[G][f]$ an $M[G][f]$-generic filter $\Gtail\of \Ptail$ and
lift the embedding to $j:V[G]\to M[j(G)]$, where $j(G)=G*f*\Gtail$. Since $\dot A\in M$, it follows that
$A\in M[j(G)]$. Thus, the function $f^*=f\union\{\<\kappa,A>\}$ is a condition in $j(\Q_\kappa)$. And since
this forcing is $\leqkappa$-closed, we may again diagonalize below the (master) condition $f^*$ to produce
an $M[j(G)]$-generic filter $j(f)$ and lift the embedding fully to $j:V[G][f]\to M[j(G)][j(f)]$. By
construction, $j(f)(\kappa)=A$, as desired. So the function $f$ is a $\LDmeas_\kappa$ Laver function.\QED

One can easily adapt this forcing to obtain the universal $\LDmeas$ Laver diamond by simply continuing the
iteration to class length. This adds a generic function $\ell$ whose restriction to any measurable cardinal
$\kappa$ is a $\LDmeas_\kappa$ Laver function. The technique adapts as well to the case of weak
compactness, unfoldability and so on.

Silver forcing is not forcing-equivalent to Woodin's fast function forcing (so the second technique I
introduced above is truly different), because Silver forcing admits a gap between any two stages of
forcing, whereas Woodin's fast function forcing has no gaps (see \cite{Hamkins2001:GapForcing}). In
particular, for any $\gamma$ above the least inaccessible cardinal, Silver forcing adds no new subsets to
$\gamma$ all of whose initial segments are in the ground model (called {\df fresh} sets in
\cite{Hamkins2001:GapForcing}); but fast function forcing adds many such sets, at any stage
$\gamma\in\dom(f)$ which is the next inaccessible cluster point $\dom(f)$ above any desired ordinal.

Finally, let me point out that in fact one can force $\LDmeas_\kappa$ with $\kappa$-c.c. forcing:

\Proof A Third Proof of Theorem \ref{FFFmeas}, with the $\kappa$-Chain-Condition: It is even possible to
force $\LDmeas_\kappa$ with $\kappa$-c.c.~forcing, if one knows $2^\kappa=\kappa^\plus$ initially. Suppose
that $\kappa$ is measurable and $2^\kappa=\kappa^\plus$. Let $\P$ be the $\kappa$-iteration that forces
with $\add(\gamma^\plus,1)$ at every inaccessible cardinal stage $\gamma<\kappa$, and suppose that $G\of\P$
is $V$-generic. If $j:V\to M$ is the elementary embedding by a normal ultrapower on $\kappa$ in $V$, then
$j(\P)$ factors as $\P*\Ptail$, where $\Ptail$ is $\leqkappa$-closed in $M[G]$. Thus, by the usual
diagonalization argument, one constructs a generic filter $\Gtail\of\Ptail$ in $V[G]$ and lifts the
embedding to $j:V[G]\to M[j(G)]$, where $j(G)=G*\Gtail$.

To see that $\LDmeas_\kappa$ holds in $V[G]$, let me extract an encoded Laver function $\ell$ from the
generic object $G$. For any inaccessible cardinal $\gamma<\kappa$, the generic object $G$ added a generic
set $G(\gamma)\of\gamma^\plus$ at stage $\gamma$. The first $\gamma$ many bits of $G(\gamma)$ code an
element $a_\gamma$ of $H(\gamma^\plus)^{V[G_\gamma]}$; let $\ell(\gamma)=a_\gamma$. I claim that
$\ell\from\kappa\to V_\kappa$ is a $\LDmeas_\kappa$ Laver function. To see this, fix any set $A\in
H(\kappa^\plus)^{V[G]}$ and any embedding $j:V\to M$ by a normal measure on $\kappa$ in $V$. Since $A$ has
a name in $H(\kappa^\plus)^V\of M$, it is in $H(\kappa^\plus)^{M[G]}$, and so in the diagonalization
argument above, we may work below a condition in $\Ptail$ such that the first $\kappa$ many bits of the
forcing at stage $\kappa$ code the set $A$. The lifted embedding $j:V[G]\to M[j(G)]$, therefore, has the
property that $j(\ell)(\kappa)=A$, as desired.\QED

This idea works also with the weakly compact cardinals, unfoldable cardinals and so on. For the larger
cardinals, such as $\theta$-supercompact cardinals, one wants to allow the stage $\kappa$ generic to code
much larger sets than merely elements of $H(\kappa^\plus)$. In order to accomplish this, simply let the
forcing at stage $\gamma$ add a subset to the next strong cardinal, say, or some other cardinals which is
guaranteed to be large enough for whatever coding is at hand.

\Section Indestructibility of the Laver diamond\hfill\break principles $\LD_\kappa$

In this section, I will prove that the various Laver diamond principles are absolute to a variety of mild
forcing extensions, such as small forcing or highly distributive forcing.

\Theorem. The Laver diamond principles $\LDwc_\kappa$, $\LDunf_\kappa$,
$\LD^{\hbox{\!\!\tiny$\Pi^m_n$-ind}}_\kappa$, $\LDsunf_\kappa$, $\LDram_\kappa$, $\LDmeas_\kappa$,
$\LDstr_\kappa$ and $\LDsc_\kappa$ are each indestructible by small forcing, that is, by forcing of size
less than $\kappa$.

\Proof: Fix a Laver function $\ell\from\kappa\to V_\kappa$ witnessing the relevant $\LDstar_\kappa$
property in $V$, and suppose that $g\of\P$ is $V$-generic for small forcing $\P\in V_\kappa$. Define
$\ell^*(\alpha)=\ell(\alpha)_g$, provided that this makes sense, that is, that $\ell(\alpha)$ is a
$\P$-name. We claim that $\ell^*$ witnesses in each case the $\LDstar_\kappa$ property in $V[g]$.

To see this, fix any appropriate $X$ that is to be guessed (e.g. for $\LDmeas_\kappa$, we should take $X\in
H(\kappa^\plus)^{V[g]}$). Let $\dot X\in V$ be a name for $X$ such that $\dot X_g=X$ and $\dot X$ has the
smallest possible hereditary size for such a name. Because the forcing is small, the set $\dot X$ is in
each case an appropriate set to be guessed by $\ell$. For example, if $X\in H(\kappa^\plus)$, then a name
$\dot X$ can also be found in $H(\kappa^\plus)$.

Let us now consider the cases of $\LDmeas_\kappa$, $\LDstr_\kappa$, $\LDsc_\kappa$ (the latter two are
trivial cases, because the large cardinal is preserved by small forcing and the Laver diamond principle
holds automatically for every strong or supercompact cardinal). Using these principles in $V$, fix an
appropriate embedding $j:V\to M$ with $j(\ell)(\kappa)=\dot X$. Standard arguments show how to lift the
embedding $j$ through the small forcing to obtain an embedding $j:V[g]\to M[j(g)]$ with $j(g)=g$ (one
simply defines $j(\tau_g)=j(\tau)_g$, and verifies that this remains elementary). Furthermore, it is well
known that this lift remains a normal ultrapower, a strongness extender embedding, or a supercompactness
embedding in $V[g]$, respectively, if $j$ had this property in $V$. Thus, the lifted embedding has the
appropriate type in $V[g]$, and by the definition of $\ell^*$ we may observe
$j(\ell^*)(\kappa)=j(\ell)(\kappa)_g=\dot X_g=X$, as desired. So $\ell^*$ is a Laver function in $V[g]$.

We now turn to the principles below a measurable cardinal, namely, $\LDwc_\kappa$, $\LDunf_\kappa$,
$\LD^{\hbox{\!\!\tiny$\Pi^m_n$-ind}}_\kappa$, $\LDsunf_\kappa$ and $\LDram_\kappa$, which follow
essentially the same argument, with a small complication. Fix any transitive set $M$ of size $\kappa$ in
$V[g]$, with $M^\ltkappa\of M$ there and $\kappa,X,\ell^*\in M$. The set $M$ has a name $\dot M$ of size
$\kappa$, and in $V$ we may find a transitive set $\bar M$ of size $\kappa$ with $\bar M^\ltkappa\of\bar M$
in $V$ and $\kappa,\dot M, \dot X,\ell\in \bar M$. For each principle, let $j:\bar M\to \bar N$ be an
embedding of the appropriate type with $j(\ell)(\kappa)=\dot X$. The standard arguments once again allow us
to lift the embedding to $j:\bar M[g]\to \bar N[j(g)]$ with $j(g)=g$. And once again, this lifted embedding
has the appropriate type in $V[G]$. Furthermore, we again know that
$j(\ell^*)(\kappa)=j(\ell)(\kappa)_g=\dot X_g=X$. Finally, we observe that the restricted embedding
$j\restrict M:M\to j(M)$ also has the appropriate type and has $j(\ell^*)(\kappa)=X$, so $\ell^*$ witnesses
that the corresponding Laver diamond principle holds in $V[g]$.\QED

The argument clearly generalizes to most all of the other Laver diamond principles; one only needs to be
able to lift the corresponding embeddings to any small forcing extension, a feature that all the large
cardinal embeddings exhibit.

Not only is the Laver diamond indestructible by small forcing, but small forcing cannot create new
instances of it.

\eject
\Theorem. Suppose that $V[g]$ is a forcing extension of\/ $V$ by forcing of size less than $\kappa$.
Then the Laver diamond principle $\LDstar_\kappa$ (in the case of weak compactness, unfoldability,
$\Pi^m_n$-indescribability, strong unfoldability, Ramseyness, measurability, superstrongness, strong
compactness and supercompactness, respectively) holds in the extension $V[g]$ if and only if it holds in
the ground model $V$.

\Proof: The converse implication is exactly the content of the previous theorem. For the forward direction,
suppose that $\ell\from\kappa\to V_\kappa$ is a Laver function for the relevant large cardinal notion in
$V[g]$, where $g\of\P$ is $V$-generic for forcing $\P$ with $\delta=|\P|<\kappa$. The function $\ell$ has a
$\P$-name $\dot\ell$ in the ground model, forced by $\one_\P$ to have these properties in the extension.
For each $p\in\P$ let $\ell_p(\alpha)=x$ if and only if $p\forces\dot\ell(\check\alpha)=\check x$.

I claim that the family of functions $\set{\ell_p\st p\in\P}$ witnesses the weak Laver diamond principle in
$V$, in the sense of Theorems \ref{WeakLDwc} and \ref{WeakLD}. Fix any set $A\in V$ to be guessed by the
Laver function. Since $\ell$ really is a Laver function in $V[g]$, there is an embedding $j$ in $V[g]$ with
$j(\ell)(\kappa)=A$. Since $j(\dot\ell)$ is a $j(\P)$-name for $j(\ell)$, there must be some condition in
$j(p)\in j(\P)=j\image\P$ forcing $j(\dot\ell)=\check A$. Thus, by $j$ of the definition of $\ell_p$, this
means that $j(\ell_p)(\kappa)=A$. Finally, since the restriction of $j$ to $V$ is an embedding in $V$ by
the Levy-Solovay Theorem, we may conclude that any set $A$ is anticipated by one of the functions $\ell_p$
in $V$. Thus, the weak Laver diamond principle holds in $V$, and so by Theorem \ref{WeakLD}, the full Laver
diamond principle holds there.\QED

The forward direction of the previous theorem generalizes to the case of forcing that is $\beta$-c.c. for
some $\beta<\kappa$. The reason is that if $\dot\ell$ is a $\P$-name for a Laver function in the extension,
then one can define $\bar\ell(\gamma)$, for $\gamma>\beta$, to be the set of all $x$ such that
$\boolval{\dot\ell(\check\gamma)=\check x}\not=0$. This will be a weak Laver function in the sense of
Theorem \ref{WeakLD}.

Let me now consider an opposite extreme of small forcing, namely, $\leqkappa$-distributive forcing, that
is, forcing that does not add any new $\kappa$-sequences over the ground model.

\Theorem. $\LDwc_\kappa$ is indestructible by $\leq\kappa$-distributive forcing, as is its
negation.\label{LDwcDistribIndestructible}

\Proof: I will show more, that the very same function witnesses $\LDwc_\kappa$ in both models. Let me
consider the weakly compact case first. Suppose that $\ell$ is a $\LDwc_\kappa$ Laver function, and that
$G\of\P$ is $V$-generic for the $\leq\kappa$-distributive forcing $\P$. I claim that $\ell$ remains a
$\LDwc_\kappa$ Laver function in $V[G]$. This is because, by the distributivity of $\P$, every set $A\in
H(\kappa^\plus)^{V[G]}$ and every nice structure $M$ of size $\kappa$ in $V[G]$ is in the ground model $V$.
Thus, there is an embedding $j:M\to N$ there with $j(\ell)(\kappa)=A$. Since this embedding still exists in
$V[G]$, the function $\ell$ witnesses $\LDwc_\kappa$ in $V[G]$.

Conversely, suppose that $\ell$ witnesses $\LDwc_\kappa$ in $V[G]$. It is easy to see that
$(V_\kappa)^{V[G]}=V_\kappa$, and so $\ell$ is a partial function from $\kappa$ to $V_\kappa$, a set in the
ground model. It follows by distributivity that $\ell$ is already in the ground model. Now, given any set
$A\in H(\kappa^\plus)$ and nice structure $M$ of size $\kappa$, there is an embedding $j:M\to N$ in $V[G]$
with $j(\ell)(\kappa)=A$. We may assume, by the usual factor arguments, that $|N|=\kappa$. By
distributivity, both the set $N$ and the embedding $j$ are in the ground model. So $\ell$ witnesses
$\LDwc_\kappa$ in $V$, as desired.\QED

\Theorem. The principles $\LDwc_\kappa$, $\LDunf_\kappa$, $\LD^{\hbox{\!\!\tiny$\Pi^m_n$-ind}}_\kappa$,
$\LDram_\kappa$, $\LDmeas_\kappa$ are each indestructible by $\leq\kappa$-distributive
forcing.\label{LDmeasDistribIndestructible}

\Proof: Again I will show the stronger result that the same Laver function serves as a witness in the
extension. To represent the smaller large cardinal cases, let me consider only $\LDunf_\kappa$. Suppose
that $\ell$ witnesses $\LDunf_\kappa$ in $V$. Fix any ordinal $\alpha$ and any nice structure $M$ of size
$\kappa$ in $V[G]$, a $\leq\kappa$-distributive forcing extension. Since $M$ must be in $V$, there is an
embedding $j:M\to N$ in $V$ with $j(\ell)(\kappa)=\alpha$. This embedding also serves as a witness in
$V[G]$, and so $\ell$ witnesses $\LDunf_\kappa$ in $V[G]$.

For the larger large cardinals, suppose that $\ell$ witnesses $\LDmeas_\kappa$ in $V$. If $A\in
H(\kappa^\plus)$ in $V[G]$, then by distributivity $A\in V$, so there is an embedding $j:V\to M$ with
$j(\ell)(\kappa)=A$. We may assume that $j$ is the ultrapower by a normal measure in $V$, so that
$M=\set{j(f)(\kappa)\st f:\kappa\to V\And f\in V}$. From this, it follows that $j\image G$ generates an
$M$-generic filter for $j(\P)$. This is because every open dense set $D\in M$ for $j(\P)$ has the form
$j(\vec D)(\kappa)$, where $\vec D=\<D_\alpha\st\alpha<\kappa>$ is a list of $\kappa$ many open dense
subsets of $\P$ in $V$. Let $\bar D=\intersect_\alpha D_\alpha$. By distributivity, this remains open and
dense. Further, since $\bar D\of D_\alpha$ for all $\alpha<\kappa$, it follows that $j(\bar D)\of D$. And
since $G$ meets $\bar D$, it follows that $j\image G$ meets $j(\bar D)$, and hence also $D$. Therefore
$j\image G$ is $M$-generic for $j(\P)$, as I claimed. Using this, one lifts the embedding to $j:V[G]\to
M[j(G)]$, where $j(G)$ is the filter generated by $j\image G$. Since $j(\ell)(\kappa)=A$ still, this
embedding witnesses $\LDmeas_\kappa$ in $V[G]$.\QED

One cannot expect to prove \ref{LDmeasDistribIndestructible} in the case of large cardinal properties that
are not themselves indestructible by $\leqkappa$-distributive forcing. For example, in the case of partial
strongness or superstrongness, while one can easily lift an extender embedding $j:V\to M$ through
$\leqkappa$-distributive forcing $G\of\Q$, since $j\image G$ generates a unique $M$-generic filter on
$j(\Q)$, the resulting lift embedding $j:V[G]\to M[j(G)]$ may not have $G\in M[j(G)]$, and therefore it may
not exhibit the required strength in $V[G]$. Indeed, it is possible to have a strong cardinal $\kappa$
whose strongness is destroyed by the forcing to add a Cohen subset to $\kappa^\plus$. Thus, one should not
expect $\LDthetastr_\kappa$ or $\LD_\kappa^{\hbox{\!\!\tiny superstrong}}$ to be preserved by
$\leqkappa$-distributive forcing. Similarly, results of \cite{Hamkins98:AsYouLikeIt} show how easy it is to
arrange that the $\kappa^\plus$-supercompactness of a cardinal $\kappa$ be destroyed after forcing to add a
subset of $\kappa^\plus$, or many other common (definable) forcing notions. So $\LDthetasc_\kappa$ is not
in general indestructible by $\leqkappa$-distributive forcing. Furthermore, one cannot hope to prove the
kind of converse to \ref{LDmeasDistribIndestructible} as we did for \ref{LDwcDistribIndestructible}, on
account of the following.

\Theorem. It is relatively consistent that $\LDmeas_\kappa$ holds in a $\leqkappa$-distributive forcing
extension, but not in the ground model. Indeed, the forcing $\add(\kappa^\plus,1)$ can create a new
instance of $\LDmeas_\kappa$.\label{NewLDmeas}

\Proof: I will build the model by forcing. Suppose in $V$ that $\kappa$ is measurable, that
$2^\kappa>\kappa^\plus$ and that $\ell$ is a $\LDmeas_\kappa$ Laver function.\footnote{The consistency
strength of this hypothesis is greater than the existence of a measurable cardinal (see
\cite{Kunen:Independence}).} Let $\P_{\kappa+1}$ be the $(\kappa+1)$-stage forcing iteration, which at
cardinal stages $\gamma\leq\kappa$, forces $2^\gamma=\gamma^\plus$ via $\add(\gamma^\plus,1)$, and suppose
that $G_{\kappa+1}\of\P_{\kappa+1}$ is $V$-generic. Factor the forcing at stage $\kappa$ as
$\P_\kappa*\Q_\kappa$, with $G_{\kappa+1}=G_\kappa*g$. The model $\bar V=V[G_\kappa]$, with its
$\leqkappa$-distributive forcing extension $\bar V[g]=V[G_\kappa][g]$, has the desired properties.

First, I claim that $\LDmeas_\kappa$ holds in $\bar V[g]$. To see this, suppose that $j:V\to M$ is any
ultrapower embedding by a normal measure on $\kappa$ in $V$. The forcing $j(\P_\kappa)$ factors as
$j(\P_\kappa)=\P_\kappa*\Q_\kappa*\Ptail$, and so we may form the partial extension $M[G_\kappa][g]$, a
model in which the forcing $\Ptail$ is $\leqkappa$-closed. Thus, by diagonalization in $V[G_\kappa][g]$,
using the fact that $2^\kappa=\kappa^\plus$ there, we may construct an $M[G_\kappa][g]$-generic filter
$\Gtail\of\Ptail$, and lift the embedding to $j:V[G_\kappa]\to M[j(G_\kappa)]$, where
$j(G)=G_\kappa*g*\Gtail$. Since the stage $\kappa$ forcing is $\leqkappa$-distributive, it follows as in
Theorem \ref{LDmeasDistribIndestructible} that $j\image g$ generates a $M[j(G_\kappa)]$-generic filter for
$j(\Q_\kappa)$, and so we may fully lift the embedding to $j:V[G_\kappa][g]\to M[j(G_\kappa)][j(g)]$. So
$\kappa$ is measurable in $\bar V[g]=V[G_\kappa][g]$. Furthermore, since the stage $\kappa$ forcing is
$\leqkappa$-distributive, any object $A\in H(\kappa^\plus)$ in $V[G_\kappa][g]$ has a $\P_\kappa$-name
$\dot A\in H(\kappa^\plus)^V$, and we could have arranged that $j(\ell)(\kappa)=\dot A$. Thus, the function
$\bar\ell(\alpha)=\ell(\alpha)_{G_\alpha}$ (provided that $\ell(\alpha)$ is a $\P_\alpha$-name) has the
property that $j(\bar\ell)(\kappa)=A$. Thus, $\bar\ell$ witnesses $\LDmeas_\kappa$ in $\bar V[g]=
V[G_\kappa][g]$. (One can alternatively argue that a Laver function is directly coded in $G_\kappa$, as in
the third proof of Theorem \ref{ForcingLD}, and avoid using $\LDmeas_\kappa$ in the ground model.)

Let me now argue that $\LDmeas_\kappa$ fails in $\bar V$. In fact, $\kappa$ is not even measurable in $\bar
V$, because the \GCH\ holds for all cardinals $\gamma<\kappa$ in $\bar V$, but not at $\kappa$ itself. In
summary, $\LDmeas_\kappa$ fails in $\bar V$, because $\kappa$ is not even measurable there, but it holds in
the $\leqkappa$-distributive extension $\bar V[g]$, as desired.\QED

The Laver function of the previous example was not actually added by the forcing $\Q_\kappa$, and indeed,
the $\leqkappa$-distributive forcing in question adds no functions on $\kappa$ at all; rather, what is
happening is that over $\bar V$ the forcing is adding the measures witnessing that the function is a Laver
function. Thus, we have the interesting situation where a function $\ell$ which is not a Laver function for
$\kappa$ in the ground model $\bar V$ becomes a Laver function in the forcing extension $\bar V[g]$. What
the theorem shows is that $\LDmeas_\kappa$ can be turned on by $\leqkappa$-distributive forcing, largely
because the measurability of $\kappa$ itself can be turned on by such forcing.  A more satisfying example
might be a model in which $\kappa$ is measurable and $\LDmeas_\kappa$ fails, but $\LDmeas_\kappa$ holds in
a $\leqkappa$-distributive forcing extension.

\Section Failures of the Laver diamond principle $\LD_\kappa$

Since the previous sections establish that the Laver diamond principle $\LDstar_\kappa$ holds in a great
variety of cases---it is nearly always forceable and often an outright consequence of the large cardinal
property $\star$ itself or a slightly stronger property---one might be tempted to expect the Laver diamond
principle for every large cardinal. In particular, one might expect that it should be difficult to make the
principle fail. But this expectation should be tempered somewhat in light of the following result, which is
especially striking in comparison with the results of Section \ref{LDinL}.

\Theorem. In the canonical inner model $L[\mu]$ for a measurable cardinal $\kappa$, the Laver diamond
$\LDmeas_\kappa$ fails.\label{LDmeasFailsInL[mu]}

\Proof: Consider the inner model $L[\mu]$, in which $\mu$ is a normal measure on $\kappa$. It is well known
that every elementary embedding definable in this model is an iteration of $\mu$. Since all such embeddings
act coherently on subsets of $\kappa$, simply stretching them taller and taller, there is for each function
$\ell\from\kappa\to V_\kappa$ only one possible value in $L[\mu]$ for $j(\ell)(\kappa)$, namely,
$j_\mu(\ell)(\kappa)$. So no function anticipates more than one set.\QED

This theorem could alternatively be deduced from the following theorem, since it is well-known that there
is only one normal measure on $\kappa$ in $L[\mu]$.

\Theorem. If there are fewer than $2^\kappa$ many normal measures on $\kappa$, then $\LDmeas_\kappa$
fails.\label{FewMeasures}

\Proof: The point is that if $\LDmeas_\kappa$ holds and $\ell$ is a $\LDmeas_\kappa$ function, then
$j(\ell)(\kappa)$ can take on any value in $H(\kappa^\plus)$, which has size $2^\kappa$. Since each of
these values requires a different normal measure, there must be at least $2^\kappa$ many normal
measures.\QED

Since all the methods we used in Section \ref{ForcingLD} to force the existence of a $\LDmeas_\kappa$ Laver
function produce models having $2^{2^\kappa}$ many normal measures on $\kappa$, one is left to wonder:

\Question. Does $\LDmeas_\kappa$ imply that there are $2^{2^\kappa}$ many normal measures on $\kappa$, the
maximal conceivable number?

Probably not, and such an answer might arise from a suitable inner model theory. Indeed, after Theorem
\ref{LDmeasFailsInL[mu]}, one wonders which of the Laver diamond principles hold in the various canonical
inner models for large cardinals. Of course, in the models where one has unique measures or extenders, then
the corresponding Laver diamond principle will fail; but what of the others?

\Question. Which of the canonical inner models for the various kinds of large cardinals satisfy the
corresponding Laver diamond principle?

Apart from the inner models, one would like to know what is possible to obtain by forcing.

\Question. If\/ $\kappa$ is a measurable cardinal, can one force $\LDmeas_\kappa$ to fail while preserving
the measurability of $\kappa$?\label{LDmeasFailQuestion}

If proper class forcing is allowed, then an affirmative answer to the question is provided by Friedman's
theorem \cite{Friedman89:CodingOverAMeasurable} that if $\kappa$ is measurable, then there is a class
forcing extension preserving this which satisfies $V=L[\mu][x]$ where $x\of\omega$. In any such model,
there is at most one normal measure on $\kappa$, and so $\LDmeas_\kappa$ will fail there. One really hopes
for a less destructive method of forcing, however, for an answer to this question, in particular, one which
has a chance to preserve large cardinals above $\kappa$. Question \ref{LDmeasFailQuestion} generalizes to
the following question:

\Question. If $\kappa$ is a $\theta$-supercompact cardinal, can one force $\neg\LDthetasc_\kappa$ while
preserving the $\theta$-supercompactness of $\kappa$? Indeed, when $\kappa<\theta$, is it relatively
consistent at all that $\kappa$ is $\theta$-supercompact but $\LDthetasc_\kappa$
fails?\label{LDthetascFailQuestion}

After Question \ref{LDthetascQuestion}, I conjectured that this should be possible. One could accomplish
this by solving another long outstanding problem: is it relatively consistent that $\kappa$ is
$\theta$-supercompact, $\kappa<\theta$, but there are relatively few (for example, fewer than $2^\theta$)
normal fine measures on $P_\kappa\theta$? If so, then $\LDthetasc_\kappa$ would fail, for the same reason
as in Theorem \ref{FewMeasures} above.

Next, I would like to point out that it is relatively consistent that $\kappa$ is weakly compact and yet
$\LDwc_\kappa$ fails.

\Theorem. If\/ $\kappa$ is weakly compact, then there is a forcing extension preserving this where
$\LDwc_\kappa$ fails. Indeed, the same result holds in the case of $\Pi^m_n$-indescribable
cardinals.\label{LDwcCanFail}

\Proof: This is an immediate consequence of Kai Hauser's \cite{Hauser92:IndescribablesWithoutDiamond}
result that one can force $\neg\Diamond_\kappa(Reg)$, while preserving the weak compactness or even the
$\Pi^m_n$-indescribability of $\kappa$.\QED

In a forthcoming article, Mirna Dzamonja and I have extended Hauser's result to the case of strongly
unfoldable cardinals, meaning that the Laver diamond $\LDsunf_\kappa$ can fail for such cardinals as well.

Let me close with the most urgent (or perhaps merely the most embarrassing) question of all.

\Question. Is $\LDwc_\kappa$ equivalent to $\Diamond_\kappa(\Reg)$ when $\kappa$ is weakly compact? Or for
that matter, is it equivalent to $\Diamond_\kappa$?

Of course, I conjecture that they are not equivalent, but no proof is forthcoming. Theorem
\ref{LDmeasFailsInL[mu]} separates $\LDmeas_\kappa$ from $\Diamond_\kappa(\Reg)$ in the case of measurable
cardinals, because it shows how $\LDmeas_\kappa$ can fail when $\kappa$ is measurable, whereas
$\Diamond_\kappa(\Reg)$ must hold whenever $\kappa$ is measurable. But we have no such result for the
smaller large cardinals, including weakly compact, unfoldable, $\Pi^m_n$-indescribable, strongly unfoldable
and Ramsey cardinals. At the moment, the only method of forcing failures of the Laver diamond principle is
to force $\neg\Diamond_\kappa(\Reg)$. Separating the notions will be a topic for further research.

\bigskip\bigskip
\noindent
jdh@hamkins.org\\
http://jdh.hamkins.org

\medskip
\noindent
{\small\sc Georgia State University\\
Department of Mathematics and Statistics\\
Atlanta, GA 30303

\medskip
{\it\ \ \&}

\medskip\noindent
The CUNY Graduate Center\\
Mathematics Program\\
365 Fifth Avenue, New York, NY 10016}

\bibliographystyle{alpha}
\bibliography{MathBiblio,HamkinsBiblio}

\begin{thebibliography}{Ham01b}

\bibitem[Cor99]{CorazzaJSL}
Paul Corazza.
\newblock Laver sequences for extendible and super-almost-huge cardinals.
\newblock {\em Journal of Symbolic Logic}, 64(3):263--283, 1999.

\bibitem[Fri89]{Friedman89:CodingOverAMeasurable}
Sy~D. Friedman.
\newblock Coding over a measurable cardinal.
\newblock {\em J. Symbolic Logic}, 54(4):1145--1159, 1989.

\bibitem[GS89]{GitikShelah89}
Moti Gitik and Saharon Shelah.
\newblock On certain indestructibility of strong cardinals and a question of
  {Hajnal}.
\newblock {\em Archive for Mathematical Logic}, 28(1):35--42, 1989.

\bibitem[Ham98]{Hamkins98:AsYouLikeIt}
Joel~David Hamkins.
\newblock Destruction or preservation as you like it.
\newblock {\em Ann. Pure Appl. Logic}, 91(2-3):191--229, 1998.

\bibitem[Ham00]{Hamkins2000:LotteryPreparation}
Joel~David Hamkins.
\newblock The lottery preparation.
\newblock {\em Ann. Pure Appl. Logic}, 101(2-3):103--146, 2000.

\bibitem[Ham01a]{Hamkins2001:GapForcing}
Joel~David Hamkins.
\newblock Gap forcing.
\newblock {\em Israel Journal of Mathematics}, 125:237--252, 2001.

\bibitem[Ham01b]{Hamkins2001:UnfoldableCardinals}
Joel~David Hamkins.
\newblock Unfoldable cardinals and the {GCH}.
\newblock {\em Journal of Symbolic Logic}, 66(3):1186--1198, 2001.

\bibitem[Hau91]{Hauser1991:IndescribableCardinals}
Kai Hauser.
\newblock Indescribable cardinals and elementary embeddings.
\newblock {\em Journal of Symbolic Logic}, 56:439--457, 1991.

\bibitem[Hau92]{Hauser92:IndescribablesWithoutDiamond}
Kai Hauser.
\newblock Indescribable cardinals without diamonds.
\newblock {\em Arch. Math. Logic}, 31(5):373--383, 1992.

\bibitem[Jec78]{Jech:SetTheory}
Thomas Jech.
\newblock {\em Set Theory}.
\newblock Academic Press, 1978.

\bibitem[Kan97]{Kanamori:TheHigherInfinite}
Akihiro Kanamori.
\newblock {\em The Higher Infinite}.
\newblock Springer-Verlag, 1997.

\bibitem[KM]{KimchiMagidor}
Y.~Kimchi and M.~Magidor.
\newblock The independence between the concepts of compactness and
  supercompactness.
\newblock Circulated manuscript.

\bibitem[Kun80]{Kunen:Independence}
K.~Kunen.
\newblock {\em Set Theory, An Introduction to Independence Proofs}.
\newblock North-Holland, 1980.

\bibitem[Lav78]{Laver78}
Richard Laver.
\newblock Making the supercompactness of $\kappa$ indestructible under
  $\kappa$-directed closed forcing.
\newblock {\em Israel Journal of Mathematics}, 29:385--388, 1978.

\bibitem[Mit79]{Mitchell1979:HypermeasurableCardinals}
William Mitchell.
\newblock Hypermeasurable cardinals.
\newblock In {\em Logic Colloquium '78 (Mons, 1978)}, volume~97 of {\em Stud.
  Logic Foundations Math.}, pages 303--316. North-Holland, Amsterdam, 1979.

\bibitem[SV]{ShelahVaananan726:ExtensionsOfInfinitaryLogic}
Saharon Shelah and Jouko V\"a\"an\"anen.
\newblock {A Note on Extensions of Infinitary Logic}.
\newblock Shelah [ShVa:726].

\bibitem[Wel02]{Welch2002:PersonalCommunication}
Philip Welch, January 2002.
\newblock personal communication.

\end{thebibliography}

\end{document}